\documentclass[11pt] {article} 
\usepackage[backref]{hyperref}
  \usepackage{amsmath}
    \usepackage{amssymb}
  \usepackage[dvipsnames]{xcolor}   
     \usepackage{pdfsync}

\newtheorem{example}{Example}[section]
\newtheorem{theorem}{Theorem}[section]
\newtheorem{lemma}{Lemma}[section]
\newtheorem{corollary}{Corollary}[section]

\newtheorem{remark}{Remark}[section]

\newcommand{\eqnsection}{
   \renewcommand{\theequation}{\thesection.\arabic{equation}}
   \makeatletter
   \csname @addtoreset\endcsname{equation}{section} 
   \makeatother}

\def \ov{\overline}

\def \be{\begin{equation}}
\def \ee{\end{equation}}
\def \bt{\begin{theorem}} 
\def \et{\end{theorem}}
\def \bl{\begin{lemma}} 
\def \el{\end{lemma}}
\def \bea{\begin{eqnarray}}
\def \eea{\end{eqnarray}}
\def \bas{\begin{eqnarray*}}
\def \eas{\end{eqnarray*}}

\def \al{\alpha}
\def \bb{\beta}
\def \ga{\gamma}

\def \de{\delta}
\def \De{\Delta}
\def \ep{\epsilon}

\def \la{\lambda}

\def \om{\omega}
\def \Om{\Omega}

\def \vf{\varphi}
\def \si{\sigma}

\def \th{\theta}

\def \ze{\zeta}

\def \ff{\infty}

\def \wt{\widetilde}

\def\stl{\stackrel{law}{=}}

\def \CC{{\cal C}}

\def \FF{{\cal F}}
\def \GG{{\cal G}}

\def \PP{{\cal P}}

\def \YY{{\cal Y}}

\def\b1{\mathbf 1}
\def \({\left(}
\def \){\right)}

\def \da{\downarrow }

\def \nn{\nonumber}
 
\def \Proof{\noindent{\bf Proof $\,$ }}

\def \bc{\begin{center} }
\def \ec{\end{center} }
\def \bs{\begin{slide} }
\def \es{\end{slide} }

\def\square{{\vcenter{\vbox{\hrule height.3pt
        \hbox{\vrule width.3pt height5pt \kern5pt
           \vrule width.3pt}
        \hrule height.3pt}}}}
\def\qed{{\hfill $\square$ \bigskip}}

\eqnsection

\begin{document}

\title{Moduli of continuity for the local times of rebirthed Markov processes }

 \author{   P.J. Fitzsimmons\,\,\, Michael B. Marcus\,\, \,\, Jay Rosen \thanks{Research of     Jay Rosen was partially supported by  grants from the Simons Foundation.   }}
\maketitle

\begin{abstract} 
Let $S$ a be locally compact space with a countable base.  
Let   $\YY $ be a transient symmetric Borel right process with state space $S$  and  continuous strictly positive  $p$--potential   densities $u^p(x,y)$.  Local and uniform moduli of continuity are obtained for the local times of both fully and partially rebirthed versions of $\YY$.  A fully rebirthed version of $\YY$ is an extension of $\YY$  so that instead of  terminating at the end of its lifetime it is immediately ``reborn'' with  a probability measure $\mu $, on $S$.   I.e., the process goes to the set $B\subset S$ with probability $\mu(B) $, after which  it continues to evolve the way $\YY$ did, being reborn with probability   $\mu $ each time it dies. 
This rebirthed version of $\YY$ is a recurrent Borel right process with state space $S$ and $p$-potential densities of form,
\[ u^p(x,y)+h(x,y),\qquad x,y\in S,\,\, p>0,
\]
where $h(x,y)$ is not symmetric.

The local times of the rebirthed process are given in terms of the local times of $\YY$  and isomorphism theorems in the spirit of Dynkin, Eisenbam and Kaspi are obtained that relate these local times to generalized chi--square processes formed by Gaussian processes with covariances $u^{q}(x,y)$ for different values of $q$. These isomorphisms allow one to obtain exact local and uniform moduli of continuity for the local times of the rebirthed  process.
Several explicit examples are given in which   $\YY$ is either a modified L\'evy process  or a  diffusion. 

Analogous results are obtained for  partially rebirthed versions of $\YY$.
 This  is  obtained by starting  $\YY$ in $S$ and  when it dies returning it to $S$  with a sub-probability measure $\Xi$.  (With probability    $1-|\Xi |$  it is sent to a disjoint state space $S'$, where it remains.)

\end{abstract}

\section{Introduction  }

  Let $S$ a be locally compact space with a countable base.  
Let   $\YY=(\Om,  \FF_{t},  \YY_t,\th_{t}, \newline P^x)$ be a transient symmetric Borel right process with state space $S$  and  continuous strictly positive  $p-$potential   densities  
\be
  u^{p}=\{ u^{p}(x,y), x,y\in S \},\quad p\geq 0,
\ee
 with respect to some $\si$--finite  positive measure $m$ on $S$. Let $\ze=\inf\{t\,|\,\YY_t=\De \}$,  where $\De$ is the cemetery state for   $\YY$  and assume       that $\ze<\ff$   almost surely.     It follows from  \cite[Lemma 3.3.3]{book} that $ \{ u^{p}(x,y), x,y\in S \} $   is positive definite and therefore is  the covariance  of a Gaussian process.
 
Let $\{G_{x}, x\in S \}$  be a mean zero Gaussian process with covariance 
 \be
 \{ u^{0}(x,y), x,y\in S \}.  
 \ee
  Isomorphism Theorems    of Dynkin and of Eisenbaum \cite[Section 8.1]{book} relate $\{G_{x}, x\in S \}$ and the   local time $\{L_{t}^{x}, x\in S$ \}   of $\YY$.  These  theorems allow us to obtain results about $\{L_{t}^{x}, x\in S \}$ from well known results about $\{G_{x}, x\in S \}$. In particular we can obtain exact local and uniform moduli of continuity for the local times of many symmetric  Borel right processes. (See   \cite[Section 9]{book}.)

  In this paper we   obtain exact local and uniform moduli of continuity for the local times of fully and partially 
 rebirthed Markov processes.   These processes do not have symmetric potential densities so the results in \cite[Section 9]{book} can not be used.
  Let $\mu$ be a probability measure on $S$. Following Meyer \cite{Meyer}        modify $\YY$ so that instead of going to $\De$ at the end of its lifetime it is immediately ``reborn'' with  probability $\mu $.  (I.e., the process goes to the set $B\subset S$ with probability $\mu(B) $, after which  it continues to evolve the way $\YY$ did, being reborn with probability   $\mu $ each time it dies.)   We denote this fully rebirthed process by $\wt Z\!=\!
(\wt \Om, \wt  \FF_{t}, \wt  Z_t, \wt \th_{t},  \wt P^y)$.

In \cite{MRLIL} we use the 
Eisenbaum--Kaspi  Isomorphism Theorem, which relates the local times of $\wt Z$ to permanental processes, to obtain exact local moduli of continuity for the local times of $\wt Z$ from  exact local moduli of continuity of associated  permanental process. 
Getting the exact local moduli of continuity for the associated  permanental process is   difficult  and the methods developed do not extend to finding uniform moduli of continuity of  permanental processes.

 Here,   using a completely different  approach, we  use properties of Gaussian processes to obtain   exact local and uniform moduli of continuity for the local times of many fully   rebirthed Markov processes $\wt  Z$.   
We show in    \cite[Theorem 7.1]{MRLIL}
that the process  
 $\wt   Z$ is a recurrent Borel right process with state space $ S$ and $p$--potential densities, 
\begin{equation}
w^{p}(x,y)=  u^{p}(x,y)+\(\frac{1}{p}-\int_S   u^{p}(x,z)\,dm(z)\)\frac{f(y)}{\|f\|_1}\label{int.1},
\end{equation}
with respect to   $m$, where  
\begin{equation} \label{int.2}
  f(y)=\int_S   u^{p}(x,y)\,d\mu(x),
\end{equation}
and the $L_1$ norm is taken with respect to $m$.
 Note that,
 \begin{equation} \label{mar.1}
  f(y)\le u^{p}(y,y),\qquad \forall y\in S,
\end{equation}  because for symmetric processes,  $u^{p}(x,y)\leq u^{p}(y,y)$,
  $\forall y\in S $. (See \cite[Lemma 3.4.3]{book}, and use $\mu(S)=1$.)

 Here is an outline of our approach.   
We have $\wt \Om=\Om^{ \mathbf N}$ with elements $\wt\om=( \om_{1},\om_{2}, \ldots)$
 and $\wt P^{y}=  P^{y}\times_{n=2}^{\ff}  P_n^{\mu}$ on $( \Om,   \FF  )^{ \mathbf N}$, where $P_n^{\mu}$ are independent copies of $P^{\mu}$,   and  where for any  Borel set $B\in S$,
 \begin{equation} \label{mar.2}
  P^{\mu}(B):=\int P^{y}(B)\,d\mu(y).
\end{equation}  
(The expectation $E^{\mu}(B)$ is defined similarly.)
  Let 
 $\ze_{n}=\ze_{n}(\wt\om)=\sum_{j=1}^{n}\ze (\om_{j})$ and set $\ze_{0}=0$.

Let $L^{y}_{t}  $  be the local time for $\YY$ normalized so that
  \begin{equation}
  E^{x}\( L_{\ff}^{y}\)= u^{0}(x,y),
  \end{equation}
and $\wt L_{t}^{y}$ denote the local time of $\wt Z$
  normalized so that,
\begin{equation}
  \wt E^{ x}\(\int_{0}^{\ff}  e^{-ps}\, d_{s}\wt L^{y}_s\)=   w^p(x,y). \label{int.3}
\end{equation}
   We show in Lemma \ref{lem-ltdecomp} that
for each $r=1,2,\ldots,$  
  \begin{equation}
 \wt L^{y}_{t}(\wt\om)=\sum^{r-1}_{i=1}L^{y}_{\ze (\om_{i})}(\om_{i})+ L^{y}_{t-\ze_{r-1}}(\om_{r}), \quad \forall t\in(\ze_{r-1},\ze_{r}], 
 \quad a.s.\label{int.4}
 \end{equation} 
This shows  that the local time $ \wt L^{y}_{t}(\wt\om)$ of $\wt Z$ can be written as a sum of independent local times $L_{\cdot}^{y}$ of the original symmetric transient process $\YY$. 

Let $\la $ denote an exponential random variable with mean $1/p$
 independent of everything else. For any    probability $P$ we denote 
 $P_{\la}=P\times \,d\la$.
 Let $\ov L^{x}_{t}$
be the local time of the process $\ov \YY=(\ov \Om, \ov  \FF_{t}, \ov  \YY_t,\ov \th_{t}, \ov  P^x)$ obtained by killing  $\YY$ at an independent  exponential random variable of mean $1/p>0$ and normalized so that
  \begin{equation}
 \ov  E^{x}\(\ov  L_{\ff}^{y}\)= u^{p}(x,y).
  \end{equation} 
 Let 
 $\{L^{x}_{i, \ff},i=1,\ldots, r-1 \}$ be independent copies of $L^{x}_{\ff}$ that are independent 
 of $\ov L^{x}_{\ff}$. 
 
 \medskip Our main result for fully rebirthed processes, given in Section \ref{sec-fullreb}, relates the local time for the rebirthed  process to the local times of the process during each cycle of rebirth.  

\medskip\noindent{\bf Theorem 2.1} {\it     For each integer $r\ge 2$ and each  $y\in S$, if   
 \begin{equation}
\(P^{y}\times \prod^{r-1}_{i=2} P^{\mu}_{i}\times \ov P^{\mu}\)  \( \sum^{r-1}_{i=1} L^{x}_{i, \ff}+\ov L^{x}_{\ff}\in B\)=1,  \label{int.6}
 \end{equation}
 for some measurable set of functions $B$,       then 
  \begin{equation}
\wt P_{\la}^{y} \(\wt L^{x}_{\la}(\wt\om)\in B\,\Big |\,    \ze_{r-1}<\la\leq \ze_{r}\)=1.\label{int.7c}
 \end{equation}
 }

As an immediate consequence of Theorem \ref{theo-1.1} and the Eisenbaum  Isomorphism  theorem,  \cite[Theorem 8.1]{book}, we  obtain the following Theorem which 
 relates Gaussian processes to the local times of fully rebirthed processes.

 \medskip\noindent{\bf Theorem 2.2} {\it For $r\ge 2$ let $\{\eta_{i, 0} (x), i  =1,\ldots,r-1 \}$  be independent  Gaussian processes with covariance $ u^{0}(x,y)$ and  $\eta_{p} (x)$ be a Gaussian process with covariance $  u^{p}(x,y)$, independent of   $\{\eta_{i, 0} (x), i  =1,\ldots,r-1  \}$. Set
\be G_{r, s}\(x\)=\sum_{i=1}^{r-1} \frac{1}{2}(\eta_{i, 0 }(x)+s)^2+\frac{1}{2}(\eta_{p}(x)+s)^2, \label{2.21pa}
\ee
and $\PP_{r,s} =  \times_{i=1}^{r-1}  P_{\eta_{i, 0} } \times   P_{\eta_{p}} $. If
\begin{equation}
\PP_{ {r, s}} \(G_{r, s}\in B\)=1,  \label{70.6nab}
 \end{equation}
  for a measurable set of functions $B\in \mathcal{M}( F(S))$, then for all $y\in S$
  \begin{equation}
\(\wt P_{\la}^{y}\times P_{G_{r, s}}\) \(\wt L^{x}_{\la}(\wt\om)+G_{r, s}\in B\,\Big |\,    \ze_{r-1}<\la\leq \ze_{r}\)=1,\qquad\forall\,x\in S.\label{70.7na}
 \end{equation}}
 
Using techniques developed in  \cite[Section 9]{book} we show that in certain circumstances 
(\ref{70.7na}) implies that,
 \begin{equation}
 \wt P_{\la}^{y}  \(\wt L^{x}_{\la}(\wt\om) \in B\,\Big |\,    \ze_{r-1}<\la\leq \ze_{r}\)=1.\label{70.7m}\ee
 This enables us to obtain properties of the local times of the fully rebirthed process. 

\medskip The condition in (\ref{70.6nab}) seems formidable. We show in Section \ref{sec-app} that it can be realized when the Gaussian processes involved satisfy relatively simple conditions. Here is an example of an application of Theorem  \ref{theo-2.1m} below  that gives an exact uniform modulus of continuity  of $\wt L^y_t$.

\medskip
Let $\eta_{i}=\{\eta_{i}(t),t\in [0,1]\}$,   $i=1,\ldots,k$, be independent   mean zero   continuous Gaussian processes with covariances $\{U_{i}(s,t),s,t\in [ 0,1]\}$ satisfying   
  \begin{equation} \label{2.38a}
 \lim_{x\to 0}\sup_{1\le i\le k}\sup_{|u-v|\le x}E(\eta_{i}(u)-\eta_{i}(v))^2\log 1/ |u-v|=0.
\end{equation}
  Let $\vf(u,v)$ be   such that for all   intervals $\De$ in $[0,1]$,
\begin{equation} \label{3.39ab}
  \vf(u,v)\le  \wt\vf\(C|u-v|\),
\end{equation}
for some constant $C$, where $\wt\vf$ is an increasing continuous function with $\wt\vf(0)=0$.

\medskip\noindent{\bf Theorem 4.2}  {\it   Assume that $S\subseteq R^{1}$.  Let $\eta_{0} (x)$ be a Gaussian process with covariance $ u^{0}(x,y)$ and let $\eta_{p} (x)$ be a Gaussian process with covariance $ u^{p}(x,y)$ and assume that (\ref{2.38a}) holds for both $\eta_0(x)$ and $\eta_p(x)$ and also that (\ref{3.39ab}) holds.       Then if,
 \be
 \lim_{h\to 0}\sup_{\stackrel{|u-v|\le h }{ u,v\in\De}} \frac{ ( \eta _{0 }(u)-\eta _{0}(v)) a_{0} + (\eta _{p }(u)-\eta _{p}(v)) a_{p} }{    \    \vf (u,v) } =1 \label{rev.1wqi} ,
\ee
  for all   intervals $\De$ in $[0,1]\cap S$ and all     $ a_0 $ and $ a_p $ with  $a_0^2+a_p^2=1$,      
    \be
\lim_{h\to 0}\sup_{\stackrel{|u-v|\le h }{ u,v\in\De}} \frac{|\wt L_{t}^{u}-\wt L_{t}^{v}|}{\vf (u,v)  }= \sup_{u\in \De}  \(2\wt L_{t}^{u}\)^{1/2},  \quad \hbox{a.e.}\,\,\, t\in S,  \,\,\, \wt P^{y}\,\,
\hbox{a.s.},\label{urev.2}
\ee
for all $y\in S$.  
 }

\medskip  In Section \ref{sec-examp} we give many examples of  Markov processes $\YY$ with potential densities that are the covariances of Gaussian processes and functions $\vf$ that satisfy (\ref{3.39ab}). Here is one of them:
   
\medskip   Let $  Y=\{Y_{t},t\in R^1_+ \}$ be a real valued symmetric 
L\'evy process with 0 potential density,   
\begin{equation} \label{1.21nk}
  \ov u(x,y)=u^\bb (x-y)=\frac{1}{2\pi}\int_{-\ff}^{\ff}\frac{ \cos\la (x-y)}{\bb+ \psi(\la)}\,d\la,\qquad x,y\in R^{1},  
\end{equation}
where $\psi$ is symmetric.    Set,
\begin{equation} \label{1.21nk}
  (\si^0)^2(x-y)=\frac{2}{ \pi}\int_{-\ff}^{\ff}\frac{1- \cos\la (x-y)}{ \psi(\la)}\,d\la,\qquad x,y\in R^{1},  
\end{equation}
  For $p>0$ let $\{\wt L_{ t}^{x},x\in R^1,t\in R^1_+ \}$ denote the local time of a process that is a full rebirthing of $Y$ determined by a probability measure $\mu$,   normalized so that,
\begin{equation}
  \wt E^{ x}\(\int_{0}^{\ff}  e^{-ps}\, d_{s}\wt L^{y}_{s}\)= u^{\bb+p}(x-y)+ \frac{\bb}{p} \int_{ -\ff} ^{\ff}  u^{\bb+p}(x,y)\,d\mu(x) ,\quad x,y\in R^1. \label{5.12}
\end{equation}

 The following result is contained in  Theorem \ref{theo-5.1}:

 \medskip\noindent{\bf Theorem $\bf (5.1)'$}   {\it Assume that $\psi$ is regularly varying at infinity with index $1<p<2$. Then        for any interval $\De$ in $[0,1]$,
   \be
\lim_{h\to 0}\sup_{\stackrel{|u-v|\le h }{ u,v\in\De}} \frac{|\wt L_{ t}^{u}-\wt L_{ t}^{v}|}{  (2(\si^0)^2(u-v) \log  1/| u-v|)^{1/2} }= \sup_{u\in \De}\(2\wt L_{ t}^{u}\)^{1/2},   \, \hbox{a.e.}\,\,\, t\in S,   \, \wt P^{y}\,\,
\hbox{a.s.}\label{urev.2}
\ee
for all $y\in R^1$.  
 }

\medskip    In Section \ref{sec-partreb} we show 
that with some minor obvious modifications all the results obtained for fully rebirthed transient Markov processes also apply to   partially rebirthed transient Markov processes.    We  define such a process as  a Borel right process
 $\wt X\!=\!
(\wt \Om,  \wt \FF_{t},\wt  X_t, \wt \th_{t},\wt P^x
)$ with state space $T=S\cup S' $,   where $ S\cap S'=\emptyset $.   We start the process in   $S$, where it behaves like $\YY $. At each death time the process is returned to $S$ with a fixed sub--probability measure $\Xi$ defined on  $S$. With probability $1-|\Xi(S)|$ it goes to $S'$, where it remains.    (Consider  $S'$ to be a place of permanent exile from   $S$.) 
We show  that it  has potential densities with a particularly simple form when restricted to $S$:
\begin{equation}
  \wt    u^{0}(x,y)= u^{0}(x,y)+f(y), \hspace{.2 in}x,y\in S,\label{int.9}
\end{equation}
where
 \begin{equation}
 f(y)=\int_S u^{0}\(x,y\)\, d\nu(x),\label{int.8}
 \end{equation}
and $\nu(x)=\Xi(x)/(1-|\Xi(S)|)$.

\section{Fully rebirthed Markov processes}\label{sec-fullreb}

    Let $S$ a be locally compact space with a countable base.  
Let   $\YY=(\Om,  \FF_{t},  \YY_t,\th_{t}, \newline P^x)$ be a transient symmetric Borel right process with state space $S$  and  continuous strictly positive  $p-$potential   densities  
\be
  u^{p}=\{ u^{p}(x,y), x,y\in S \},\quad p\geq 0,
\ee
 with respect to some $\si$--finite  positive measure $m$ on $S$. Let $\ze=\inf\{t\,|\,\YY_t=\De \}$,  where $\De$ is the cemetery state for   $\YY$  and assume       that $\ze<\ff$   almost surely.    
  Let  $L_t^y=\{L^{y}_{t},y\in S,t\in R^+\}  $ be the local time for $\YY$ normalized so that
  \begin{equation}\label{2.2}
  E^{x}\( L_{\ff}^{y}\)= u^{0}(x,y).
  \end{equation}
  We require that $L^{y}_{t}  $ is jointly continuous; see Remark \ref{rem-3.1}.
 
 \medskip Let $\mu$ be a probability measure on $S$. Following \cite{Meyer} we       modify $\YY$ so that instead of going to $\De$ at the end of its lifetime it is immediately ``reborn'' with  measure $\mu $.  (I.e., the process goes to the set $B\subset S$ with probability $\mu(B) $, after which  it continues to evolve the way $\YY$ did, being reborn with probability   $\mu $ each time it dies.)  
  
   We denote this rebirthed process by $\wt Z\!=\!
(\wt \Om, \wt  \FF_{t}, \wt  Z_t, \wt \th_{t},\wt P^x)$.   
We have $\wt \Om=\Om^{ \mathbf N}$ with elements $\wt\om=( \om_{1},\om_{2}, \ldots)$
 and $\wt P^{y}=  P_{1}^{y}\times_{n=2}^{\ff}P_n^{\mu}$ on $( \Om,   \FF  )^{ \mathbf N}$, where   for any probability measure $P$, $\{P_i,i\ge 1 \}$ are independent copies of $P$. 
  Let 
 $\ze_{n}=\ze_{n}(\wt\om)=\sum_{j=1}^{n}\ze (\om_{j})$ and set $\ze_{0}=0$.

We show in    \cite[Theorem 7.1]{MRLIL}
that the process  
 $\wt   Z$ is a recurrent Borel right process with state space $ S$ and $p$--potential densities, 
\begin{equation}
w^{p}(x,y)=  u^{p}(x,y)+\(\frac{1}{p}-\int_S   u^{p}(x,z)\,dm(z)\)\frac{f(y)}{\|f\|_1}\label{193.3a},
\end{equation}
with respect to   $m$, where  
\begin{equation} \label{7.9a}
  f(y)=\int_S   u^{p}(x,y)\,d\mu(x),
\end{equation}
and the $L_1$ norm is taken with respect to $m$.
  We point out in    (\ref{mar.1})   that $f(y)<\ff$ for all $y\in S$.

 \bl\label{lem-ltdecomp}  The process $\wt Z$ has a jointly continuous local time $\{\wt L_{t}^{y}, y\in S,t\in R_+ \}$ such that
for each $r\ge 1$,
  \begin{equation}
 \wt L^{y}_{t}(\wt\om)=\sum^{r-1}_{i=1}L^{y}_{\ze (\om_{i})}(\om_{i})+ L^{y}_{t-\ze_{r-1}}(\om_{r}), \quad \forall t\in(\ze_{r-1},\ze_{r}], 
 \quad a.s.\label{2.12t}
 \end{equation} 
 and for each $p>0$
 \begin{equation}
  \wt E^{ x}\(\int_{S}   e^{-ps}\, d_{s}\wt L^{y}_s\)=   w^p(x,y). \label{az.1vw}
\end{equation}
  \el
 
 \Proof   Since $L_{0}^{y}(\om_{r})=0$, (\ref{2.12t}) is equivalent to showing that for each $r=1,2,\ldots$
  \begin{equation}
 \wt L^{y}_{t}(\wt\om)=\sum^{r-1}_{i=1}L^{y}_{\ze (\om_{i})}(\om_{i})+ L^{y}_{t-\ze_{r-1}}(\om_{r}), \quad \forall t\in[\ze_{r-1},\ze_{r}), 
 \quad a.s.\label{2.12tr}
 \end{equation}
    Let   \begin{equation}
N_t:= N_{t}(  \wt\om)=\min \{j\,|\,t< \ze_{j} (\wt \om)\}. \label{}
  \end{equation}
  We define the stochastic process  $\{A^{y}_{t},y,t\in S\times R^1_+ \}$ as follows:
  For $\ze_{r-1}\leq t<\ze_{r}$, so that   $N_{t}=r$,  
  \begin{equation}
 A^{y}_{t}:=A^{y}_{t}(\wt\om)= \sum^{ N_{t}-1}_{i=1}L^{y}_{\ze (\om_{i})}(\om_{i})+ L^{y}_{t-\ze_{ N_{t}-1}}(\om_{ N_{t}}).\label{rcaf.1}
  \end{equation} 
Thus, for $\ze_{r-1}\leq  t< \ze_{r}$ 
 \begin{equation}
 A^{y}_{t}(\wt\om)= \sum^{ r-1}_{i=1}L^{y}_{\ze (\om_{i})}(\om_{i})+ L^{y}_{t-\ze_{ r-1}}(\om_{ r}).\label{caffy.1}
 \end{equation}
 Note that $ A^{y}_{t}$ is a continuous function of $t$. For $\ze_{r-1}< t<\ze_{r}$ this follows from the continuity of $L^{y}_{t}$. The continuity at $ t=\ze_{r-1}$ follows because $L_{0}^{y}(\om_{r})=0$ so that 
 \be A^{y}_{\ze_{r-1}}(\wt\om)= \sum^{ r-1}_{i=1}L^{y}_{\ze (\om_{i})}(\om_{i}),
 \ee 
 while, by the analogue of (\ref{caffy.1}) for $\ze_{r-2}\leq  t< \ze_{r-1}$
  \begin{equation}
 A^{y}_{t}(\wt\om)= \sum^{ r-2}_{i=1}L^{y}_{\ze (\om_{i})}(\om_{i})+ L^{y}_{t-\ze_{ r-2}}(\om_{ r-1}),\label{caffy.1b}
 \end{equation}
 and $\lim_{t\to \ze_{ r-1}} L^{y}_{t-\ze_{ r-2}}(\om_{ r-1})=L^{y}_{\ze (\om_{r-1})}(\om_{r-1})$.

\medskip  We now show that  $\{A^{y}_{t}, t\geq 0 \}$ is the local time of $\wt Z$. 
    We first show that $A^{y}_{t}$ is an additive functional for $\wt   Z$, that is,
\begin{equation}
A^{y}_{t+s}=A^{y}_{s}+A^{y}_{t}\circ \wt \theta_{s}.\label{rcaf.10}
\end{equation}
By (\ref{rcaf.1}) this is equivalent to showing that
\begin{eqnarray}
&&\sum^{ N_{t+s}-1}_{i=1}L^{y}_{\ze (\om_{i})}(\om_{i})+ L^{y}_{t+s-\ze_{ N_{t+s}-1}}(\om_{ N_{t+s}})
\label{rcaf.50}
\\
&&\qquad =\sum^{ N_{s}-1}_{i=1}L^{y}_{\ze (\om_{i})}(\om_{i})+ L^{y}_{s-\ze_{ N_{s}-1}}(\om_{ N_{s}})+A^{y}_{t}\circ \wt \theta_{s}.
\nonumber
\end{eqnarray}
We show below that 
\be A^{y}_{t}\circ \wt \theta_{s}
\label{rcaf.13r}
 =L^{y}_{\ze (\om_{N_s})}(\om_{N_s})\circ  \theta_{s-\ze_{N_s-1}}+\sum^{ N_{s+t}-1}_{i=N_s+1}L^{y}_{\ze (\om_{i})}(\om_{i}) +L^{y}_{s+t- \ze_{N_{s+t}-1} }(\om_{N_{s+t}}).
\ee
Then since
\begin{equation}
L^{y}_{\ze (\om_{N_s})}(\om_{N_s})=L^{y}_{s-\ze_{ N_{s}-1}}(\om_{ N_{s}})+L^{y}_{\ze (\om_{N_s})}(\om_{N_s})\circ  \theta_{s-\ze_{N_s-1}}\label{},
\end{equation}
we see that (\ref{rcaf.13r}) implies (\ref{rcaf.50}).

\medskip We now prove (\ref{rcaf.13r}).
Suppose that $N_s=m$, which implies that that $\ze_{m-1}\leq s< \ze_{m}$. 
Then
\be
\wt \theta_{s}( \om_{1},\om_{2}, \ldots)=(  \theta_{s-\ze_{m-1}}\om_{m},\om_{m+1}, \om_{m+2},\ldots \,), \label{rcaf.11}
\ee
so that 
\begin{equation}
\ze_{j}\circ \wt \theta_{s}=\ze_{j+m-1}-s.\label{rcaf.51}
\end{equation}
 To see this note that the first death time on the right of (\ref{rcaf.11}), \be \ze_{1}\circ \wt \theta_{s}=\ze(\om_m)-(s-\ze_{m-1})=\ze_{m}-s
. \ee
The next death time comes by adding $\ze(\om_{m+1})$ to this to obtain $\ze_{m+1}-s$. Continuing this process gives (\ref{rcaf.51}).
(Observations similar to (\ref{rcaf.51}) and (\ref{rcaf.11}) are given in    \cite[(5)]{Meyer} and in the proof of the Lemma in \cite[page 469 ]{Meyer},   with a   different notation.)  
Therefore 
  \bea
 N_{t}\circ \wt \theta_{s}&=&\min \{j\,|\,t< \ze_{j}\circ \wt \theta_{s} \} \label{rcaf.52}\\
&=&\min \{j\,|\,s+t< \ze_{j+m-1}  \}= N_{s+t}-\(m-1\).\nn
  \eea
Using (\ref{rcaf.11}) we now see that,
  \bea
&&\wt \theta_{s}( \om_{1},\om_{2}, \ldots,\om_{N_{t}-1}, \om_{N_{t}} )\label{rcaf.53}\\
&&\qquad=(  \theta_{s-\ze_{m-1}}\om_{m},\om_{m+1}, \om_{m+2},\ldots,  \om_{N_{t}\circ \wt \theta_{s}-1+m-1}, \om_{N_{t}\circ \wt \theta_{s}+m-1}) \nn\\
&&\qquad=(  \theta_{s-\ze_{N_s-1}}\om_{N_s},\om_{N_s+1}, \om_{N_s+2},\ldots, \om_{N_{s+t}-1},\om_{N_{s+t}}).\nn
\eea
 In particular,
\begin{equation}
\wt \theta_{s}\om_{1}=  \theta_{s-\ze_{N_s-1}}\om_{N_s},\label{rcaf.53a}
\end{equation}
\begin{equation}
\wt \theta_{s}\om_{2}=\om_{N_s+1},\,\,\,\,\ldots \ldots\,\,\,\, ,\wt \theta_{s}\om_{N_{t}-1}=\om_{N_{s+t}-1},\wt \theta_{s}\om_{N_{t}}=\om_{N_{s+t}}\label{rcaf.53b}
\end{equation}
and
\begin{equation}
\ze_{ N_{t}-1}\circ \wt \theta_{s}=\ze_{\(N_{t}\circ \wt \theta_{s}-1\)}-s=\ze_{N_{s+t} -1}-s.\label{rcaf.57}
\end{equation}

By (\ref{rcaf.1})
  \bea
 A^{y}_{t}\circ \wt \theta_{s}&=& \sum^{ N_{t}-1}_{i=1}L^{y}_{\ze (\om_{i})}(\om_{i})\circ \wt \theta_{s}+ L^{y}_{t-\ze_{ N_{t}-1}}(\om_{ N_{t}})\circ \wt \theta_{s} \\
&=&L^{y}_{\ze (\om_{1})}(\om_{1})\circ \wt \theta_{s}+ \sum^{ N_{t}-1}_{i=2}L^{y}_{\ze (\om_{i})}(\om_{i})\circ \wt \theta_{s}+ L^{y}_{t-\ze_{ N_{t}-1}}(\om_{ N_{t}})\circ \wt \theta_{s}.\nn
 \eea
Using this and (\ref{rcaf.53a})--(\ref{rcaf.57})  we obtain (\ref{rcaf.13r}).
 
\medskip Therefore $A^{y}_{t}$ is a continuous additive functional for $\wt   Z$.  Let $R_{A^{y}}=\inf \{t\,|\, A^{y}_{t}>0\}$. It is clear that $\wt P^{y}\(R_{A^{y}}=0\)=1 $ and $\wt P^{x}\(R_{A^{y}}=0\)=0$ for all $x\neq y$.  Therefore  $A^{y}_{t}$ is a local time for $\wt   Z$.   (See,   e.g., \cite[(3.82)]{book}).

	\medskip  We now show that $A^{y}_{t}$ has the normalization given in (\ref{az.1vw}).      By (\ref{caffy.1}), for $\ze_{r-1}\leq  t< \ze_{r}$, 
 \begin{equation}
 A^{y}_{t}(\wt\om)= \sum^{ r-1}_{i=1}L^{y}_{\ze (\om_{i})}(\om_{i})+ L^{y}_{t-\ze_{ r-1}}(\om_{ r}).\label{caffy.1a}
 \end{equation}
Therefore, for $\ze_{r-1}< t< \ze_{r}$,
 \begin{equation}
 d_{t}A^{y}_{t}=d_{t}L^{y}_{t-\ze_{r-1}}(\om_{r}).\label{caffy.2}
 \end{equation}
%

 Since  $ A^{y}_{t}$ is continuous in $t$,  $ d_{t}A^{y}_{t}$  puts no mass at the points    $\ze_{1}, \ze_{2}, \ldots$.        Consequently,
  \begin{eqnarray}
  &&  \wt E^{ x}\(\int_{0}^{\ff}  e^{-pt}\,  d_{t}A^{y}_{t}\)=\sum_{r=1}^{\ff}  \wt E^{ x}\(\int_{(\ze_{r-1},\ze_{r})}  e^{-pt}\,  d_{t}L^{y}_{t-\ze_{r-1}}(\om_{r})\)
  \label{caf.3}
  \\
  &&\qquad=\sum_{r=1}^{\ff}  \wt E^{ x}\( e^{-p\ze_{r-1}}    \int_{0}^{\ze (\om_{r})}  e^{-pt}\,  d_{t}L^{y}_{t}(\om_{r})\)
  \nonumber \\
  &&\qquad= E^{x} \(  \int_{0}^{\ff}  e^{-pt}\,  d_{t}L^{y}_{t}\)+\sum_{r=2}^{\ff}  \wt E^{ x}\( e^{-p\ze_{r-1}}\) E^{ \mu} \(  \int_{0}^{\ff}  e^{-pt}\,  d_{t}L^{y}_{t}\).
  \nonumber \eea 
 Using  \cite[Theorem 3.6.3]{book} we see that  this 
  \bea 
  &&\qquad=  u^{p}(x,y)+\sum_{r=2}^{\ff}  \wt E^{ x}\( e^{-p\ze_{r-1}}\)\int_S u^{p}(x,y)\,d\mu(x)
  \\
  &&\qquad=  u^{p}(x,y)+\sum_{r=2}^{\ff}  \wt E^{ x}\( e^{-p\ze_{r-1}}\)f(y). 
  \nonumber
  \end{eqnarray}
   To show that this gives (\ref{az.1vw}) it remains to show that, 
  \begin{equation}
 \sum_{r=2}^{\ff}  \wt E^{ x}\( e^{-p\ze_{r-1}}\)= \(\frac{1}{p}-\int_S   u^{p}(x,z)\,dm(z)\)\frac{1}{\|f\|_1}.\label{rcaf.12s}
  \end{equation}
   We have 
  \begin{equation}
 \sum_{r=2}^{\ff}  \wt E^{ x}\( e^{-p\ze_{r-1}}\)=\sum_{r=1}^{\ff}  \wt E^{ x}\( e^{-p\ze_{r}}\)=\sum_{r=1}^{\ff}  \wt E^{ x}\( e^{-p\sum_{j=1}^{r}\ze (\om_{j})}\). \label{}
  \end{equation}
Using the Markov property it is easy to see that 
  \begin{equation}
  \wt E^{ x}\( e^{-p\sum_{j=1}^{r}\ze (\om_{j})}\)=E^{ x}\( e^{-p \ze}\)\(E^{\mu}\( e^{-p \ze}\)\)^{r-1}\label{}
  \end{equation}
so that 
\begin{equation}
 \sum_{r=2}^{\ff}  \wt E^{ x}\( e^{-p\ze_{r-1}}\)=\frac{E^{ x}\( e^{-p \ze}\)}{1-E^{\mu}\( e^{-p \ze}\)}.\label{}
\end{equation}
By simple integration,
  \begin{equation}
E^{x}\( e^{-p\ze}\)=1-pE^{x}\( \int_{0}^{\ze}e^{-pt} \,dt \)=1-p\int_{S} u^{p}(x,z)\,dm(z),  \label{rp.4}
\end{equation}
where we use,  
\begin{equation}
E^{x}\( \int_{0}^{\ze}e^{-pt} \,dt\)=E^{x}\( \int_{0}^{\ff}e^{-pt}1_{\{\YY_t\in S\}} \,dt\)
=\int_{0}^{\ff}\int_{S}e^{-pt}p_{t}(x,z)  \,dm(z) \,dt.\label{}
\end{equation}
  Now it follows from (\ref{rp.4}) that,
\begin{equation} 
  E^{\mu}\( e^{-p\ze}\)= 1-p\int_{S} \int_{S} u^{p}(x,z)\,dm(z)\,d\mu(x)=1-p\|f\|_{1},  \label{rp.4a}
\end{equation}
so that $1-E^{\mu}\( e^{-p \ze}\)=p\|f\|_{1}$.
This shows that (\ref{rcaf.12s}) holds.  

 Set  $\wt  L^{x}_{t}(\wt\om)=A^{y}_{t}(\wt\om)$. To complete the proof  it only remains to
show  that $\wt  L^{x}_{t}$ is jointly continuous. By (\ref{2.12t}),
   \begin{equation}
 \wt L^{x}_{t}(\wt\om)=\sum^{r-1}_{i=1}  L^{x}_{\ze(\om_{i})}(\om_{i})+   L^{x}_{t-\ze_{r-1}}(\om_{r}),\qquad \ze_{r-1}<t\leq \ze_{r}.\label{70.1jc}
 \end{equation}
 Consequently, the joint continuity follows immediately from the  joint continuity of  $ L^{x}_{t}$ except when  $t=\ze_{r}$. By (\ref{70.1jc})
    \begin{equation}
 \wt L^{x}_{\ze_{r}}(\wt\om)=\sum^{r}_{i=1}  L^{x}_{\ze(\om_{i})}(\om_{i}),\label{70.1jd}
 \end{equation}
 and, 
    \begin{equation}
 \wt L^{x}_{\ze_{r}+s}(\wt\om)=\sum^{r}_{i=1}  L^{x}_{\ze(\om_{i})}(\om_{i})+   L^{x}_{s}(\om_{r+1}),\qquad 0<s\leq \ze(\om_{r+1}).\label{70.1je}
 \end{equation}
Since $L^{x}_{0}\equiv 0$ this shows that we have joint continuity at $t=\ze_{r}$.\qed

By (\ref{2.12t}), when $\ze_{r-1}<t\leq \ze_{r}$, and $r\geq 1$, 
 we have
  \begin{equation}
 \wt L^{x}_{t}(\wt\om)=\sum^{r-1}_{i=1}  L^{x}_{\ze(\om_{i})}(\om_{i})+   L^{x}_{t-\ze_{r-1}}(\om_{r}).\label{70.1}
 \end{equation} 
  We now show that the terms on the right in (\ref{70.1}) are conditionally independent. Let $\la,  \la_{1}, \la_{2},\ldots  $ denote  independent exponential random variables of mean $1/p$
 independent of everything else. For a probability measure $P$   we write 
 $P_{\la}=P\times \,d\la$ and similarly, 
 $P_{\la_{i}}=P\times \,d\la_{i}$.
 Let  $\la' =\la-\ze_{r-1}$, when $\ze_{r-1}<\la\leq \ze_{r}$,
   \begin{equation}
 \wt L^{x}_{\la}(\wt\om)=\sum^{r-1}_{i=1}  L^{x}_{\ze(\om_{i})}(\om_{i}) +     L^{x}_{\la'\wedge \ze(\om_{r})}(\om_{r}).\label{70.2}
 \end{equation}
  Let
 \bea
 \mathcal{L}_{i}&=&\{L^{x}_{\ze(\om_{i})}(\om_{i}), x\in S\},\qquad1\le i\leq r-1,  \label{}\\
 \mathcal{L}'_{r}&=&\nn \{L^{x}_{\la'\wedge \ze(\om_{r})}(\om_{r}), x\in S\},\label{}\\
   \ov {\mathcal{L}} _{r} &=&\nn \{L^{x}_{\la_{r} }(\om_{r}), x\in S\}.\label{}\eea
In this notation (\ref{70.2}) can be written as,
     \begin{equation}
 \wt L^{x}_{\la}(\wt\om)=\sum^{r-1}_{i=1}   \mathcal{L}_{i} +  \mathcal{L}'_{r}.\label{70.2a}
 \end{equation}

\bl\label{lem-condind}
 For all $r\geq 2$ and all Borel sets, $B_1, B_2,\ldots, B_r$, in $C(S, R^1)$,
  \begin{eqnarray}
 &&\wt P_{\la}^{y}\( \mathcal{L}_{1}\in B_1, \ldots,  \mathcal{L}_{r-1}\in B_{r-1},  \mathcal{L}'_{r} \in B_r \,\Big |\,    \ze_{r-1}<\la\leq \ze_{r} \)
\nn
 \\
 &&\qquad= P_{1,\la_{1}}^{y}\( \mathcal{L}_{1}\in B_1 \,\Big |\,   \ze(\om_{1})<\la_{1}   \) \prod^{r-1}_{i=2}   P_{i,\la_{i}}^{\mu}\(\mathcal{L}_{i}\in B_i\,\Big |\,  \ze(\om_{i})<\la_{i} \)\nn\\
 &&\qquad\qquad   \times      P_{r,\la_{r}}^{\mu}\( \ov {\mathcal{L}}_{r}\in B_r\,\Big |\,  \la_{r}\leq \ze (\om_{r})\).
 \label{70.5mp}
 \end{eqnarray}  
 In particular,
 \begin{equation} \label{3.26}
   \wt P_{\la}^{y}\(  \mathcal{L} '_{r}\in B_r\,\Big |\, \ze_{r-1}  <\la\le  \ze _{r}\)=  P_{r,\la_{r}}^{\mu}\( \ov {\mathcal{L}}_{r}\in B_r\,\Big |\,  \la_{r}\leq \ze (\om_{r})\),
\end{equation}
which implies that for each $r\geq 2$,  
 \begin{equation}
 \mathcal{L}_{1}, \,\ldots, \, \mathcal{L}_{r-1},\,  \mathcal{L}'_{r}, \label{70.0}
 \end{equation}
 are conditionally independent, given that $\ze_{r-1}<\la\leq \ze_{r}$. 

  \el
  
 \medskip
   \Proof Let $\mathcal{C}$ be a countable subset of $S$ with compact closure.    For any function $f(x)$,   $x\in\CC$, and $\nu_i\in R^1$, $1\le i\le r$,   set,  
\begin{equation} \label{}
   (\nu_{i},f)=\sum_{x\in \mathcal{C}}a_{i, x} f(x),
\end{equation}  where the $a_{i, x}\in R^{1}$, $1\le i\le r$    and $\sum_{x\in \mathcal{C}}|a_{i, x}|<\ff$  for each $i$. 
We begin this proof by calculating the Laplace transform of,  
\begin{equation} \label{}
  \-\sum^{r-1}_{i=1} (\nu_{i},   L^{\cdot}_{\ze(\om_{i})}(\om_{i}))+  (\nu_{r},  L^{\cdot}_{\la'\wedge \ze(\om_{r})}(\om_{r})),
\end{equation} 
where $L^\cdot$ is equal $L^x$ for some $x\in S.$
Using (\ref{70.2})  we have,
 \begin{eqnarray}
 \lefteqn{\wt E_{\la}^{y}\(\exp\({-\sum^{r-1}_{i=1} (\nu_{i},   L^{\cdot}_{\ze(\om_{i})}(\om_{i}))-  (\nu_{r},  L^{\cdot}_{\la'\wedge \ze(\om_{r})}(\om_{r}))}\);    \ze_{r-1}<\la\leq \ze_{r} \)}
 \label{70.3}
 \\
 &&=\wt E^{y} \int_{ \ze_{r-1}}^{ \ze_{r-1}+\ze (\om_{r})} pe^{-pt} \exp\({-\sum^{r-1}_{i=1} (\nu_{i},  L^{\cdot}_{\ze(\om_{i})}(\om_{i}))-  (\nu_{r},  L^{\cdot}_{t-  \ze_{r-1}}(\om_{r}))}     \) \,dt 
 \nonumber \\
 &&=\wt E^{y}\int_{ 0}^{ \ze (\om_{r})} pe^{-p(u+ \ze_{r-1})}  \exp \(-\sum^{r-1}_{i=1} (\nu_{i},  L^{\cdot}_{\ze(\om_{i})(\om_{i}))}-  (\nu_{r},   L^{\cdot}_{u}(\om_{r}))    \)  \,du.\nn    
 \eea
  
Since, obviously, for any function $f(u)$
\begin{equation} \label{}
  \int_{ 0}^{ \ze (\om_{r})} pe^{-p(u+ \ze_{r-1})}f(u) \,du=e^{-p \ze_{r-1} } \int_{ 0}^{ \ze (\om_{r})} pe^{-p u }f(u) \,du,
\end{equation}
we can write the last line of (\ref{70.3}) in the form,
\be 
 \wt E^{y}\( e^{-p \ze_{r-1}} e^{-\sum^{r-1}_{i=1} (\nu_{i},  L^{\cdot}_{\ze(\om_{i})}(\om_{i}))}    \int_{ 0}^{ \ze (\om_{r})} pe^{-pu} e^{-  (\nu_{r},  L^{\cdot}_{u}(\om_{r}))}   \,du\).
  \ee
Using the Markov property this is
\begin{equation} \label{}
  \wt E^{y}\( \prod^{r-1}_{i=1}  e^{-p \ze(\om_{i}) } e^{-  (\nu_{i},  L^{\cdot}_{\ze(\om_{i})}(\om_{i}))}     \)        E^{\mu}\(\int_{ 0}^{ \ze (\om_{r})} pe^{-pu} e^{- (\nu_{r},   L^{\cdot}_{u}(\om_{r}))}   \,du\).
\end{equation}
This is equal to
 \bea  &&   E_{1,\la_{1}}^{y}\(  e^{-  (\nu_{1},  L^{\cdot}_{\ze(\om_{1})}(\om_{1}))} ;  \ze(\om_{1})<\la_{1}   \) \prod^{r-1}_{i=2}   E_{i,\la_{i}}^{\mu}\( e^{-  (\nu_{i},  L^{\cdot}_{\ze(\om_{i})}(\om_{i}))} ; \ze(\om_{i})<\la_{i} \) \nonumber\\
 &&\hspace{1 in} \times      E_{r.\la_{r}}^{\mu}\( e^{-  (\nu_{r},   L^{\cdot}_{\la_{r}}(\om_{r}))}; \la_{r}\leq \ze (\om_{r})\),\label{70.3a}
 \end{eqnarray}
 where $E^{\cdot}_{i,\la_i}$ is expectation with respect to $P_{i,\la_i}$.
 
  Using (\ref{70.3}) and (\ref{70.3a})  and setting all the $\nu_{i}=0$ we see that,
 \be \wt P_{\la}^{y} \(  \ze_{r-1}<\la\leq \ze_{r}\)=   P_{1,\la_{1}}^{y} \( \ze(\om_{1})<\la_{1} \)\prod^{r-1}_{i=2}   P_{i,\la_{i}}^{\mu} \( \ze(\om_{i})<\la_{i}\)   P_{r,\la_{r}}^{\mu} \(  \la_{r}\leq \ze (\om_{r})\),
 \ee
 Consequently,
 \begin{eqnarray}
 \lefteqn{\wt E_{\la}^{y}\(e^{-\sum^{r-1}_{i=1} (\nu_{i},  L^{\cdot}_{\ze(\om_{i})}(\om_{i}))-  (\nu_{r},  L^{\cdot}_{\la'\wedge \ze(\om_{r})}(\om_{r}))} \,\Big |\,    \ze_{r-1}<\la\leq \ze_{r} \)
 \label{70.5}}
 \\
 &&=  E_{1,\la_{1}}^{y}\(  e^{-  (\nu_{1},  L^{\cdot}_{\ze(\om_{1})}(\om_{1}))} \,\Big |\, \ze(\om_{1})<\la_{1}   \) \prod^{r-1}_{i=2}   E_{i,\la_{i}}^{\mu}\( e^{-  (\nu_{i},  L^{\cdot}_{\ze(\om_{i})}(\om_{i}))}  \,\Big |\, \ze(\om_{i})<\la_{i} \) \nonumber\\
 &&\hspace{1 in} \times      E_{r,\la_{r}}^{\mu}\( e^{-  (\nu_{r},   L^{\cdot}_{\la_{r}}(\om_{r}))} \,\Big |\, \la_{r}\leq \ze (\om_{r})\).
 \nonumber
 \end{eqnarray}
 
 Taking the inverse  Laplace transform    gives  (\ref{70.5mp}) with the $B_{i}$  restricted to $\mathcal{C}$. The statement in (\ref{70.0}) follows because the local times are continuous. \qed

 Let 
 $ L^{x}_{i, t}$    be independent copies of $L^{x}_{t}$  and  let $\ze(i)$ denote the death time of  $L^{x}_{i,t}$, $i=1, 2, \ldots  $. 
 For any set $C$, let $F(C)$ denote the set of real--valued functions
$f$ on $C$. Define the evaluations
$\it{i}_{ x}:F(C)\mapsto R^{{\rm 1}}$ by $\it{i}_{ x}( f)=f( x)$. We use
$\mathcal{M}( F(C))$ to denote  the smallest $\si$-algebra for which the
evaluations
$\it{i}_{ x}$ are Borel measurable for all
$x\in C$. ($\mathcal{M}( F(C))$ is generally referred to  as the
$\si$-algebra  of cylinder sets in $F(C)$.)

 \bt\label{theo-1.1} When, for some $r\geq 2$, and measurable set of functions $B\in \mathcal{M}( F(S))$,
 \begin{equation}
\(P_{1}^{y}\times \prod^{r-1}_{i=2} P_{i}^{\mu}\times   P_{r, \la_{r}}^{\mu}\)  \( \sum^{r-1}_{i=1} L^{\cdot}_{i, \ff}+  L^{\cdot}_{r,\la_{r}}\in B\)=1,  \label{70.6}
 \end{equation}
 then  
  \begin{equation}
\wt P_{\la}^{y} \(\wt L^{\cdot}_{\la}(\wt\om)\in B\,\Big |\,    \ze_{r-1}<\la\leq \ze_{r}\)=1.\label{70.7}
 \end{equation}
 \et
 
 \Proof  
 Since $\{\la_{i}, 1\le i\le r-1\}$ do not appear in the event  $\{ \sum^{r-1}_{i=1} L^{\cdot}_{i, \ff}+  L^{\cdot}_{r,\la_{r}}\in B\}$, we can  write  (\ref{70.6}) as,
  \begin{equation}
\(P_{1,\la_{1}}^{y}\times \prod^{r-1}_{i=2} P_{i,\la_{i}}^{\mu}\times   P_{r, \la_{r}}^{\mu}\)  \( \sum^{r-1}_{i=1} L^{\cdot}_{i, \ff}+  L^{\cdot}_{r,\la_{r}}\in B\)=1. \label{70.6c}
 \end{equation}
Now using the fact that  an almost sure event   with respect to a given probability occurs almost surely for any conditional version of the probability, we see 
  that   (\ref{70.6c}) implies that,     
  \be\(P_{1,\la_{1}}^{y}\times \prod^{r-1}_{i=2} P_{i,\la_{i}}^{\mu}\times   P_{r, \la_{r}}^{\mu}\) \( \sum^{r-1}_{i=1} L^{\cdot}_{i, \ff}+  L^{\cdot}_{r,\la_{r}}\in B\,\Big |\,\mathcal{A}\)=1,
  \label{70.6d}
\ee 
 where
  \begin{equation}
 \mathcal{A}=\{ \ze(1)<\la_{1},\,\cdots,\,  \ze(r-1)<\la_{r-1},  \,  \la_{r}\leq \ze (r)\}.    \label{}
  \end{equation}
 We note that   for each $i$, $\mathcal{L}_{i}= L^{x}_{\ze(\om_{i})}(\om_{i})$ has the law of $L^{x}_{\ff}$ 
 and $ \ov {\mathcal{L}}= L^{x}_{\la}(\om_{r})$  has the law of $ L^{x}_{\la}$. Hence (\ref{70.6d}) is equivalent to 
    \begin{eqnarray}
  &&\(P_{1,\la_{1}}^{y}\times \prod^{r-1}_{i=2} P_{i,\la_{i}}^{\mu}\times   P_{r, \la_{r}}^{\mu}\)
  \label{70.6e}
  \\
  &&\qquad  \(   \sum^{r-1}_{i=1} \mathcal{L}_{i}+\ov {\mathcal{L}}\in B\,\Big |\,    \ze(\om_{1})<\la_{1},\cdots,  \ze(\om_{r-1})<\la_{r-1},  \,  \la_{r}\leq \ze (\om_{r}) \)=1.
  \nonumber
  \end{eqnarray}
It then follows from (\ref{70.5mp}) that, 
  \begin{equation}
  \wt P_{\la}^{y}\(\sum^{r-1}_{i=1} \mathcal{L}_{i}+\mathcal{L}'_{r}\in B   \,\Big |\,    \ze_{r-1}<\la\leq \ze_{r} \)=1.\label{}
  \end{equation}
Using (\ref{70.2a})  we see that this is the statement in (\ref{70.7}).
\qed

\bt\label{theo-ITcond}  Let $\{\eta_{i, 0} (x)$, i=1,\ldots\} be independent  Gaussian processes with covariance $ u^{0}(x,y)$. Let   $\eta_{p} (x)$ be a Gaussian process with covariance $  u^{p}(x,y)$ independent of the $\{\eta_{i, 0} (x)\}$. For any $r\geq 1$ set $G_{r, s}=\{G_{r, s}(x),x\in S \}$ where
\be G_{r, s}\(x\)=\sum_{i=1}^{r-1} \frac{1}{2}(\eta_{i, 0 }(x)+s)^2+\frac{1}{2}(\eta_{p}(x)+s)^2, \label{2.21p}
\ee
and $ P_{G_{r, s}} =  \times_{i=1}^{r-1}  P_{\eta_{i, 0} } \times   P_{\eta_{p}} $. When 
\begin{equation}
 P_{G_{r, s}} \(G_{r, s}\in B\)=1,  \label{70.6n}
 \end{equation}
  for a measurable set of functions $B\in \mathcal{M}( F(S))$, then 
  \begin{equation}
\(\wt P_{\la}^{y}\times P_{G_{r, s}}\) \(\wt L^{x}_{\la}(\wt\om)+G_{r, s} \in B\,\Big |\,    \ze_{r-1}<\la\leq \ze_{r}\)=1,\qquad   \forall y\in S.\label{70.7n}
 \end{equation}
\et

\noindent\textbf{Proof of Theorem \ref{theo-ITcond} }
 It follows from the Eisenbaum Isomorphism Theorem, \cite[Theorem 8.1.1]{book}, that,   
\bea\lefteqn{
\Big\{ L^x_{i,\ff}+\frac{1}{2}(\eta_{i, 0} (x)+s)^2\,\,;\,x\in S\,,\,P_{i}^y\times
P_{\eta_{i, 0}}\Big\}\label{it1.2f}}\\ &&\qquad\stackrel{law}{=}
\Big\{\frac{1}{2}(\eta_{i, 0} (x)+s)^2\,\,;\,x\in S\,,\,\(1+\frac{\eta_{0,i} (y)}s\)P_{\eta_{i, 0}}\Big\},\nn
\eea
where $P_{\eta_{i, 0}}$ is the probability  of $\eta_{i, 0}$ and $P_{i}^y$ the probability of $\YY$ started at $y$. 
It also follows from the Eisenbaum Isomorphism Theorem that,
\bea\lefteqn{
\Big\{ L^x_{r, \la_{r}}+\frac{1}{2}(\eta_{p} (x)+s)^2\,\,;\,x\in S\,,\,P_{r, \la_{r}}^y\times
P_{\eta_{p}}\Big\}\label{it1.2}}\\ &&\qquad\stackrel{law}{=}
\Big\{\frac{1}{2}(\eta_{p} (x)+s)^2\,\,;\,x\in S\,,\,(1+\frac{\eta_{p } (y)}s)P_{\eta_{p}}\Big\}\nn,
\eea
where $P_{\eta_{p}}$ is the probability of $\eta_{p}$  and $P_{r, \la_{r}}^y= P_{r}^y\times d\la_{r}$. 

Set,
\begin{equation} \label{3.37} 
  \eta_{i, 0} (\mu)=\int_S \eta_{i, 0} (x)\,d\mu(x),\qquad\text{and}\qquad \eta_{p} (\mu)=\int_S \eta_{p} (x)\,d\mu(x).
\end{equation}   It follows (\ref{it1.2f}) and (\ref{it1.2})   that  for each  $r\geq 2$,
\bea
&& 
\Big\{ \sum^{r-1}_{i=1} L^{x}_{i, \ff}+ L^{x}_{r, \la_{r}}+\sum_{i=1}^{r-1}\frac{1}{2}(\eta_{i, 0}(x)+s)^2+\frac{1}{2}(\eta_{p}(x)+s)^2\,;\,\, x\in S,\,\,\nn  \\
&&\hspace{1 in}  \( P_{1}^{y} \times\prod _{i=2}^{r-1}P_i^{\mu}\times P_{r, \la_{r}}^{\mu} \)     \times\prod_{i=1}^{r-1}  P_{\eta_{i, 0} } \times   P_{\eta_{p} } \Big\} \label{3.38}\\
 &&\stackrel{law}{=}
\Big\{ \sum_{i=1}^{r-1}\frac{1}{2}(\eta_{i, 0}(x)+s)^2+\frac{1}{2}(\eta_{p}(x)+s)^2\,\,\,\,\,;\,x\in S\,,\nn\\
&& \hspace{.5 in}\qquad     \,(1+\frac{\eta_{ 0,1}(y)}s)P_{\eta_{0,1} } \times\prod_{i=2}^{r-1}  (1+\frac{\eta_{i, 0}(\mu)}s)P_{\eta_{i, 0} }\times   ( (1+\frac{\eta_{p}(\mu)}s)P_{\eta_{p} }) \Big\},\nn
\eea
  where $\{P^{\mu}_{i},1\le i\le r-1 \}$  are all equal to $P^\mu$.
    Note that the statement of this theorem for  $r=1$ is precisely (\ref{it1.2}).  
  
  For $r>1$ it follows from (\ref{70.6n}) and (\ref{3.38}) that 
 \be
 \(P_{1}^{y} \times_{i=2}^{r-1}P_i^{\mu}\times  P^{\mu}_{r,\la_{r}}   \times   P_{G_{r, s} }\)
\( \sum^{r-1}_{i=1} L^{\cdot}_{i, \ff}+ L^{\cdot}_{r,\la_{r}}+G_{r,s}(\cdot)\in B\)=1. \label{}
\ee
Therefore, for $P_{G_{r, s} }$ almost every $\om'$,
 \be
 \(P_{1}^{y} \times_{i=2}^{r-1}P_i^{\mu}\times  P^{\mu}_{r,\la_{r}}  \)
\( \sum^{r-1}_{i=1} L^{\cdot}_{i, \ff}+ L^{\cdot}_{r,\la_{r}}+G_{r,s}(\cdot, \om')\in B\)=1. \label{}
\ee
 This  theorem now follows    from   Theorem \ref{theo-1.1}.\qed  
 
 \begin{remark} {\rm \label{rem-3.1} Theorem 
  \ref{theo-ITcond} relates local times of fully rebirthed Markov processes  to certain Gaussian processes. Our applications of this theorem require that the Gaussian processes have continuous sample paths. It follows from  \cite[Theorem 9.4.1]{book} that the associated local times are jointly continuous. 

 }\end{remark}

 \section{Moduli of continuity of generalized chi-square processes }\label{sec-mod}  

The expression in (\ref{2.21p}) is formidable. In this section we show how it can be simplified for  the applications of Theorem \ref{theo-ITcond} in this paper.

\subsection{Local modulus of continuity }

  Let $\eta_{i}=\{\eta_{i}(t),t\in [0,1]\}$, $i=0,\ldots,k$, be independent    continuous  mean zero Gaussian processes with covariances $U_{i}=\{U_{i}(s,t),s,t\in [0,1]\},$ with $U_{i}(d,d)>0$, where $d\in [0,1].$   Consider, 
 \begin{equation} \label{1.9mm}
  Y_{k}(t)= \sum_{i=1}^{k} \eta^2_{i}(t),\qquad t\in [0,1].\end{equation}
  The stochastic process $Y_{k } =\{Y_{k }(t),t\in  [0,1] \}$ is referred to as a     generalized   chi--square process of order $k $.  (We say generalized because in the standard definition of a  chi--square process of order $k $ the $\{\eta_i \}_{i=1}^k$     are   identicaly distributed.)

  Set
   \begin{equation}
  \si_{i}^2(s,t)=E(\eta_{i}(s)-\eta_{i}(t))^2, \qquad i=0,\ldots k.\label{II.10}
   \end{equation}
  We assume that, 
      \begin{equation}
 \lim_{s,t\to d} \frac{\si_{i}^2(s,t)}{\si_{0} ^2(s,t)}=1, \qquad i=1, \ldots,k \label{1.4a},
   \end{equation}
  and that  $\{\eta_0
  (t),t\in [0,1] \}$ has a local modulus of continuity   at $d\in [0,1]$,
  \begin{equation} \label{1.7}
   \limsup_{t\to 0}  \frac{ |\eta_0(t+d)-\eta_0 (d)|} {   \phi(t)} =  1    \qquad a.s.,
\end{equation}
for some increasing function $\phi$.

\medskip When $\{\eta_{i}(t),t\in  R^1\}$,    $i=0,\ldots,k$ are stationary Gaussian processes we sometimes write,
\begin{equation} \label{}
  U_{i}(s-t):=U_{i}(s,t)\qquad\text{and}\qquad   \si^2_{i}(s-t):=\si^2_{i}(s,t), \qquad i=1, \ldots,k.
\end{equation}
 In this case we can write (\ref{1.4a}) as,
   \begin{equation}\label{1.6}
 \lim_{ t\to 0} \frac{\si_{i}^2(t)}{\si_{0} ^2( t)}=1, \qquad i=1, \ldots,k .
 \ee

  \bt \label{theo-1.3m}   When the conditions in (\ref{1.4a}) and (\ref{1.7})  are satisfied, 
  \begin{equation} 
\limsup_{t\to 0} \frac{|Y_{k }(t+d)-Y_{k }(d)|} {   \phi(t)} =   2 Y^{1/2}_{k}(d) \qquad a.s.\label{rev.2qq}
\end{equation} 
 \et

 \Proof  To simplify the notation we initially take $d=0$.   This proof is a generalization of the proof of  \cite[Theorem 2.1]{MRECP}. We   refer to that proof for some details.
Let   $\mathbf{a}=\{a_i,i=1,\ldots,k \}$ be a sequence of real numbers with $\|\mathbf{a}\|_2=1$ and consider,
\begin{equation} \label{1.09}
  \ov\eta (t) :=  \sum_{i=1}^{k}  a_i \eta_{i} (t). 
\end{equation}
It follows from (\ref{1.4a}) that for all $\ep>0$, there exists a $\de>0$ such that for  $s,t\in [0,\de]$,
\begin{equation} \label{1.11}
(1-\ep)\si^{2}_0(s,t) \le  E(\ov \eta (s)-\ov \eta (t))^2\le (1+\ep)\si^{2}_0(s,t).
\end{equation}
It follows from this, (\ref{1.7}) and  \cite[Lemma 7.1.10 and Remark 7.1.11]{book} that,
  \begin{equation} \label{1.7dd}
\limsup_{t\to 0} \frac{|\ov \eta (t)-\ov \eta  (0)|} {   \phi(t)} =  1    \qquad a.s.
\end{equation}
(See Remark \ref{rem-1.1} for more details.)  
We write,   
\bea \label{}
 {\ov\eta (t)-\ov\eta (0) }&:= & \sum_{i=1}^{k}  a_i\(  \eta_{i} (t)-\eta_{i} (0) \)\\
&=& \sum_{i=1}^{k} a_i \(  \eta_{i} (t)-\frac{U_{i}(0,t)}{U_{i}(0,0)} \eta_{i} (0)\)   \nn\\
&& \qquad-\sum_{i=1}^{k}a_i\(\frac{U_{i}(0,0)-U_{i}(0,t)}{U_{i}(0,0)} \) \eta_{i} (0).\nn
\eea

It follows from (\ref{1.7dd}) that for each $i\in {1,\ldots, k}$
  \begin{equation} \label{1.7dd4}
\limsup_{t\to 0} \frac{|  \eta_{i} (t)-  \eta_{i}  (0)|} {   \phi(t)} =  1    \qquad a.s.
\end{equation}
Consequently, by  \cite[({2.5})]{MRECP} we see that,  
\bea \label{deltau.1}
 &&  \lim_{t\to 0} \frac{ | U_{i}(0,0)-U_{i}(0,t)|}{\phi(t)}       =0.
\eea

Therefore,
 \begin{equation} \label{}
\limsup_{t\to 0}\frac{|\sum_{i=1}^{k}a_i\(\frac{U_{i}(0,0)-U_{i}(0,t)}{U_{i}(0,0)} \) \eta_{i} (0)|}{\phi(t)}=0,
\end{equation}
and consequently,
\begin{equation} \label{1.20}
\limsup_{t\to 0}\frac{\big|\sum_{i=1}^{k} a_i \(  \eta_{i} (t)-\frac{U_{i}(0,t)}{U_{i}(0,0)} \eta_{i} (0)\)\big|}{\phi(t)}=1.
\end{equation}

 \medskip   Consider the stochastic process,
 \begin{equation} \label{}
  \rho(t):=  \sum_{i=1}^{k}\( \eta_{i}(t)-\frac{U_{i}(0,t)}{U_{i}(0,0)} \eta_{i}(0)\)  \eta_{i}(0) ,\qquad t\in [0,1].
\end{equation}

 Let $(\Om, \mathcal{F},\mathbb{P})$ be the sample space generated by the Gaussian random variables $\{\eta_{i}(t), t\in [0,1], i=1,\ldots, k\}$  and      $(\wt \Om , \wt {\mathcal{F}} ,\wt {\mathbb{P}} )$ be the sample space generated by the Gaussian random variables \[\bigg\{\eta_{i} (t)-\displaystyle \frac{U_{i}(0,t)}{U_{i}(0,0)} \eta_{i} (0), t\in [0,1], i=1,\ldots, k \bigg \}.\] Let $(R^{k}, \mathcal{E},\mu)$ be the sample space generated by the Gaussian random variables $\{\eta_{i}(0), i=1,\ldots, k \}$, Then  since 
   $\( \eta_{i}(t)-(U_{i}(0,t)/U_{i}(0,0)) \eta_{i}(0)\)$  and $\eta_{i}(0)$ are independent for $i=1,\ldots,k$,
 we have, 
  \begin{equation}
(\Om, \mathcal{F},\mathbb{P})=(\wt \Om , \wt {\mathcal{F}} ,\wt {\mathbb{P}} )\times (R^{k}, \mathcal{E},\mu). 
 \label{prod.1}
 \end{equation} 
  For each     $\vec \eta (0, \om')=(\eta_{1} (0,\om'), \dots, \eta_{k}(0,\om'))\in (R^{k}, \mathcal{E},\mu)$ consider,  
   \be 
   \rho_{\om'}(t):= \sum_{i=1}^{k}  \(  \eta_{i} (t)-\frac{U_{i}(0,t)}{U_{i}(0,0)} \eta_{i} (0)\)    \eta_{i}(0,\om').\qquad   t\in [0,1].  \label{2.10mm}
      \ee 
It now follows from (\ref{1.20})  that  
   \begin{equation} \label{1.7aa}
\wt {\mathbb{P}} \( \limsup_{t\to 0} \frac{|\rho_{\om'} (t)| } {   \phi(t)} =   \|\vec \eta(0,\om')\|_2\)=1. 
  \end{equation} 
 It follows from this, (\ref{prod.1})  and Fubini's Theorem that,  
  \begin{equation} \label{1.7w}
  \mathbb{P}  \( \limsup_{t\to 0}\frac{|\rho  (t)| } {   \phi(t)} =   \|(\eta_{1} (0), \dots, \eta_{k}(0))\|_2\)=1. \end{equation} 
Since $\|(\eta_{1} (0), \dots, \eta_{k}(0))\|_2=Y^{1/2}(0)$, we can write this as,
\begin{equation} \label{1.7x}
  \mathbb{P}  \(  \limsup_{t\to 0}\frac{|\rho  (t)| } {   \phi(t)} =   Y^{1/2}(0) \)=1. \end{equation} 
 We show in  \cite[(2.6)]{MRECP} that,  
\be 
 \limsup_{t\to 0}\frac{|Y_{k }(t)-Y_{k }(0)|}{ | \phi (t)| }\label{5.4mmt}\\
=  \limsup_{t\to 0}\frac{ 2 |\rho  (t)|}{   \phi (t)  }.
   \ee
   (Actually in \cite[(2.6)]{MRECP} the $\eta_i$ are taken to be i.i.d. It is easy to see that the independence of the $\eta_i$ suffices.)

  The statements in (\ref{1.7x}) and (\ref{5.4mmt}) give (\ref{rev.2qq}) with $d=0$. It is easy to extend it to hold for all $d\in [0,1]$. Consider,
\begin{equation} \label{}
  \eta(t+d):=\wt \eta(t),\quad t\in [0,1-d],
\end{equation}
where $\{\wt \eta(t), t\in [0,1-d] \}$ satisfies the hypotheses of this theorem when $d=0$.
We have just proved that (\ref{rev.2qq}) holds for $\wt\eta(t)$. Substituting $  \eta(t+d) $ for $\wt \eta(t)$ we get (\ref{rev.2qq}) for general $d\in R^1.$
\qed

 \begin{remark}{\rm  \label{rem-1.1}
 It follows from (\ref{1.11})  and   \cite[Lemma 7.1.10 and Remark 7.1.11]{book} that, 
 
  \begin{equation} \label{1.24}
   \limsup_{t\to 0}\sup_{ u\le t\atop u\in[0,\de]} \frac{|\ov \eta (u)-\ov \eta  (0)|} {   \phi(t)} =  1    \qquad a.s.,
\end{equation}
if and only if,
  \begin{equation} \label{1.25}
   \limsup_{t\to 0}\sup_{ u\le t\atop u\in[0,\de]} \frac{ | \eta_0 (u)- \eta_0  (0)|} {   \phi(t)} =  1    \qquad a.s.
\end{equation}
 
Furthermore, when $\phi(t)$ is increasing we now show that we can write   (\ref{1.25}) as,
  \begin{equation} \label{1.27}
   \limsup _{t\to 0}  \frac{ | \eta_0 (t)- \eta_0  (0)|} {   \phi(t)} =  1    \qquad a.s.
\end{equation}
 By   (\ref{1.25}) the left--hand side of (\ref{1.27}) is clearly less than or equal to 1. It also follows from   (\ref{1.25}) that for any sequence $\{\ep_k \}\da 0$ there exist $\{u_k \}$ and $\{t_k \} $, with $\{u_k\le t_k \}$ such that,  
\begin{equation} \label{}
   \frac{ | \eta_0 (u_k)- \eta_0  (0)|} {   \phi(t_k)}\ge 1-\ep_k.
\end{equation}
Since $\{u_k\le t_k \}$ and $\phi(t)$ is increasing this implies that, 
\begin{equation} \label{}
   \frac{ | \eta_0 (u_k)- \eta_0  (0)|} {   \phi(u_k)}\ge 1-\ep_k.
\end{equation}
This shows that the left--hand side of (\ref{1.27}) is greater  than or equal to 1. Consequently, we get (\ref{1.27}). A similar argument shows that (\ref{1.27}) implies 
(\ref{1.25}). The same analysis will work for (\ref{1.24}). 

We note that by  \cite[Theorem 7.1.4]{book} whenever, $\si_i(t,0)$ is asymptotic 
to an increasing function we can find a local modulus $\phi(t)$ that is increasing.
To simplify our expressions we have inserted   in Theorem \ref{theo-1.3m} the hypotheses that 
$\phi(t)$ is increasing. 

}\end{remark}

 The following corollary of  Theorem \ref{theo-1.3m}   is   used in our applications of the Eisenbaum Isomorphism, \cite[Theorem 8.1.1]{book} to relate generalized chi--square processes to the local times of related Markov processes. 
 
   \begin{corollary} \label{cor-2.1}   When the conditions in (\ref{1.4a}) and (\ref{1.7})  are satisfied, for any  $s\in R^{1}$, $d\in [0,1)$, 
\begin{equation} 
 \limsup_{t\to 0} \frac{\sum_{i=1}^{k}( \eta_{i}(t+d+s)^2- \sum_{i=1}^{k}( \eta_{i}(d)+s)^2} {   \phi(t)} =   2 \(\sum_{i=1}^{k}( \eta_{i}(d)+s)^2\)^{1/2} , \label{1.35ksa}
\ee  almost surely.
 \end{corollary}

 \Proof  As in the proof of Theorem \ref{theo-1.3m} we initially take $d=0$.     Let $\ov \eta_{i}=\{\eta_{i}(t)+\xi,t\in [0,1]\}$,   $i=0,\ldots,k$, where $\xi$ is a standard normal random variable independent of  $\eta_{i}=\{\eta_{i}(t),t\in [0,1]\}$,    $i=0,\ldots,k$.   Clearly,   
\begin{equation}
  \si_{i}^2(u,v):=E(\ov \eta_{i}(u)-\ov \eta_{i}(v))^2=E(\eta_{i}(u)-\eta_{i}(v))^2, \qquad i=0,\ldots k.\label{II.10q}
   \end{equation}
   Also, trivially,
   \begin{equation} \label{1.7qq}
   \limsup_{t\to 0} \frac{\ov\eta_0(t)-\ov\eta_0 (0)} {   \phi(t)} =  1    \quad a.s. \iff   \limsup_{t\to 0} \frac{ \eta_0(t)-\eta_0 (0)} {   \phi(t)} =  1    \quad a.s.
\end{equation}
Note that,
\begin{equation} \label{1.9mmq}
  Y_{k}(t):= \sum_{i=i}^{k}  \ov \eta _{i}^{2}(t)  = \sum_{i=i}^{k}( \eta _{i}(t)+\xi)^2,\qquad t\in [ 0,1].
  \ee
  It follows from Theorem \ref{theo-1.3m} that,
\begin{equation} 
 \limsup_{t\to 0} \frac{Y_{k }(t)-Y_{k }(0)} {   \phi(t)} =   2 Y^{1/2}_{k}(0) \qquad a.s.,\label{rev.2ll}
\end{equation} 
  or, equivalently,
 \begin{equation} 
 \limsup_{t\to 0} \frac{\sum_{i=i}^{k}( \eta _{i}(t)+\xi)^2- \sum_{i=i}^{k}(  \eta  _{i}(0)+\xi)^2} {   \phi(t)} =   2 \(\sum_{i=i}^{k}( \eta _{i}(0)+\xi)^2\)^{1/2}  \quad a.s.\label{rev.2ll}
\end{equation}
  It follows from this for any countable dense   set $S'$ of $R^1$, 
\begin{equation} 
 \limsup_{t\to 0} \frac{\sum_{i=i}^{k}(   \eta  _{i}(t)+s')^2- \sum_{i=i}^{k}(   \eta  _{i}(0)+s')^2} {   \phi(t)} =   2 \(\sum_{i=i}^{k}(   \eta  _{i}(0)+s')^2\)^{1/2}  \label{1.35}
\end{equation}
for all $s'\in S'$  almost surely. 

For any $s\in R^1$ and $1\le i\le k$,
\bea \label{1.36}
   && (\eta _{i}(t)+s)^2-  ( \eta _{i}(0)+s)^2 \\
   &&\qquad\nn  =(\eta _{i}(t)+s')^2-  ( \eta _{i}(0)+s')^2 -2(s'-s)(\eta_i(t)-\eta_i(0)).
\eea
Denote the left--hand side of (\ref{1.35}) by 
${\cal N}_{k,s'}$. Using (\ref{1.36}) we see that any $s\in R^1$,
\begin{equation} \label{}
  \big |{\cal  N}_{k,s}-{\cal  N}_{k,s'}\big |\le 2|s'-s|\limsup_{t\to 0} \frac{\sum_{i=i}^{k}  \eta _{i}(t) -   \eta _{i}(0) } {   \phi(t)}  \le 2k|s'-s|,  
\end{equation}
where we use (\ref{1.09})--(\ref{1.7dd}) with 
 $a_i=1  /\sqrt k $ for $1\le i\le k$.  Using this and (\ref{1.35})  we see that (\ref{1.35}) holds for all $s\in R^1$.
 
 Using the argument at the end of the proof of Theorem \ref{theo-1.3m} we can generalize this to (\ref{1.35ksa})  when $d\neq 0$.
 
  \subsection{Uniform modulus of continuity }

We now give an analogue of Theorem \ref{theo-1.3m} for the uniform modulus of continuity. This is more complicated. A major problem is that we don't have an analogue of  \cite[Lemma 7.1.10]{book} for the uniform modulus of continuity. Nevertheless, a similar result does hold for certain classes of Gaussian processes.

 Let $\eta_{i}=\{\eta_{i}(t),t\in [0,1]\}$,   $i=1,\ldots,k$, be independent   mean zero   continuous Gaussian processes with covariances $\{U_{i}(s,t),s,t\in [ 0,1]\}$ satisfying   
  \begin{equation} \label{2.38}
 \lim_{x\to 0}\sup_{1\le i\le k}\sup_{|u-v|\le x}E(\eta_{i}(u)-\eta_{i}(v))^2\log 1/ |u-v|=0.
\end{equation}Define $ \GG_{\vf,\De,k}$ as all such sequences  $\{\eta_{i} \}_{i=1}^k$  for which,
\be
\lim_{h\to 0}\sup_{\stackrel{|u-v|\le h }{ u,v\in\De}} \frac{\sum_{i=1}^{k}( \eta _{i }(u)-\eta _{i }(v)) a_{i}   }{    \    \vf (u,v) } =1 \label{rev.1wquqq} , \quad\hbox{ a.s. }
\ee
  for all   intervals $\De$ in $[0,1]$ and all     $\{a_{i} \}_{i=1}^k$  with  $\| (a_{1},\ldots , a_{k}) \|_{2}=1$.

   Let $\vf(u,v)$ be   such that 
\begin{equation} \label{3.39}
  \vf(u,v)\le  \wt\vf\(C|u-v|\),
\end{equation}
for some constant $C$, where $\wt\vf$ is an increasing  continuous function with $\wt\vf(0)=0$. 

\bt  \label{theo-2.1m}  Let  $  \eta=\{\eta_{i}(t),t\in [0,1],i=1,\ldots,k\}\in   \GG_{\vf,\De,k}$.     Then,

\begin{equation} 
\lim_{h\to 0}\sup_{\stackrel{|u-v|\le h }{ u,v\in\De}}  \frac{|Y_{k }(u)-Y_{k }(v) |}{ \vf (u,v )} =   2 \sup_{u\in\De}Y_{k }^{1/2}(u), \qquad a.s.\label{rev.2quu}
\end{equation} 
and for any $s\in R^{1}$,  
\begin{equation} 
\lim_{h\to 0}\sup_{\stackrel{|u-v|\le h }{ u,v\in\De}}  \frac{\sum_{i=i}^{k}( \eta_{i}(u)+s)^2- \sum_{i=i}^{k}( \eta_{i}(v)+s)^2} {   \vf (u,v)} =   2\sup_{u\in\De} \(\sum_{i=i}^{k}( \eta_{i}(u)+s)^2\)^{1/2} , \label{rev.2qusa}
\ee  almost surely.

\et

\Proof   Assume first that $U^*:=\inf_{i=1,\ldots,k,\,t\in[0,1]} U_{i}(t,t)>0$. 
To show that,
   \begin{equation}
\lim_{h\to 0}\sup_{\stackrel{|u-v|\le h }{ u,v\in\De}}  \frac{|Y_{k }(u)-Y_{k }(v) |}{   \vf  (u,v)} \geq   2 \sup_{t\in\De}Y_{k }^{1/2}(t), \hspace{.2 in}a.s.\label{rev.2qu5}
\end{equation} 
it suffices to show that for any $d\in \De$,
\begin{equation}
\lim_{h\to 0}\sup_{\stackrel{|u-v|\le h }{ u,v\in\De}} \frac {|Y_{k }(u)-Y_{k }(v) |}{    \vf   (u,v)} \geq  2  Y_{k }^{1/2}(d), \hspace{.2 in}a.s.\label{rev.2qu6}
\end{equation} 
This is because (\ref{rev.2qu6}) holding almost surely implies that for any countable dense set $\De'\subset \De$,
\begin{equation}
\lim_{h\to 0}\sup_{\stackrel{|u-v|\le h }{ u,v\in\De}}  \frac{|Y_{k }(u)-Y_{k }(v)| }{    \vf   (u,v)} \geq 2\sup_{d\in\De'}    Y_{k }^{1/2}(d), \hspace{.2 in}a.s.\label{rev.ioi}
\end{equation} 
which implies (\ref{rev.2qu5}).

Consider,
\be  \label{}
{\ov\eta (u)-\ov\eta (v) } :=  \sum_{i=1}^{k}  a_i\(  \eta_{i} (u)-\eta_{i} (v) \),
 \ee 
and write,
\begin{equation} \label{2.50mm}
 \eta_{i}(u)- \eta_{i}(v) = V_{i}(v,d) \eta_{i}(u)-V_{i}(u,d)\eta_{i}(v) +  G_{i}(u,v),
\end{equation}
where, $V_{i}(u,v)= {U_{i}(u,v)}/{U_{i}(d,d)}$, 
and
\begin{equation} \label{1.45}
 G_{i}(u,v)=(1-V_{i}(v,d) )\eta_{i}(u) -(1-V_{i}(u,d))\eta_{i}(v) .
\end{equation}

We show in Lemma \ref{lem-3.1} that,   
 \begin{equation}
\lim_{h\to 0}\sup_{\stackrel{|u-v|\le h }{ u,v\in\De}}  \frac{| \sum_{i=1}^k a_i G_{i}(u,v)| }{    \vf  (u,v)}\le\frac{k ^{1/2}\wt\vf( C|\De| )}{(U^*)^{1/2}},\qquad a.s. \label{1.44}
\end{equation} 
 Using this we see that,  
\begin{eqnarray}
&& \lim_{h\to 0}\sup_{\stackrel{|u-v|\le h }{ u,v\in\De}} \frac{\sum_{i=1}^{k}( \eta _{i }(u)-\eta _{i }(v)) a_{i}   }{    \   \vf (u,v)} \label{2.50nn3} \\
&&\qquad    \le \lim_{h\to 0}\sup_{\stackrel{|u-v|\le h }{ u,v\in\De}} \frac{\sum_{i=1}^{k}\(V_{i }(v,d) \eta_{i}(u)-V_{i }(u,d)\eta_{i}(v)\) a_{i} }{     \vf (u,v)}+ \frac{k ^{1/2} \wt \vf(C|\De| )}{(U^*)^{1/2}}.\nn
\end{eqnarray}
By   (\ref{rev.1wquqq}) the left--hand side of (\ref{2.50nn3}) is equal to 1. Therefore, 
\be 
\lim_{h\to 0}\sup_{\stackrel{|u-v|\le h }{ u,v\in\De}} \frac{\sum_{i=1}^{k}\(V_{i }(v,d) \eta_{i}(u)-V_{i }(u,d)\eta_{i}(v)\) a_{i} }{   \vf (u,v)}\ge 1-\frac{k ^{1/2}\wt \vf(C|\De|)}{(U^*)^{1/2}} .\label{2.50nn4} 
\ee
Note that for all $u,v\in [0,1],$
\begin{equation}  
 E\((V_{i }(v,d) \eta_{i}(u)-V_{i }(u,d)\eta_{i}(v))\eta_{i}(d)\)=0.
\end{equation}
This shows that $\eta_{i}(d)$ is independent of $\{ V_{i }(v,d) \eta_{i}(u)-V_{i }(u,d)\eta_{i}(v);u,v\in [0,1]\}$. 
It follows from this and the paragraph containing (\ref{prod.1}) that,
\bea &&\wt {\mathbb{P}}  \( \lim_{h\to 0}\sup_{\stackrel{|u-v|\le h }{ u,v\in\De}} \frac{\sum_{i=1}^{k}\(V_{i }(v,d) \eta_{i}(u)-V_{i }(u,d)\eta_{i}(v)\)\eta_{i}(d,\om') }{     \vf (u,v)    \|\vec \eta(d,\om') \|_2 }\right.\\
&&\left.\hspace{3in}
 \ge 1-\frac{k ^{1/2}\wt \vf(C|\De|) }{(U^*)^{1/2}}\)=1 \nn.\label{} 
\eea
As in  (\ref{1.7aa}) and (\ref{1.7w}) we   then  have,\bea\label{1.51} &&\mathbb{P}\( \lim_{h\to 0}\sup_{\stackrel{|u-v|\le h }{ u,v\in\De}} \frac{\sum_{i=1}^{k}\(V_{i }(v,d) \eta_{i}(u)-V_{i }(u,d)\eta_{i}(v)\) \eta_{i}(d) }{     \vf (u,v)} \right.\\
&&\hspace{1in}\left.\ge Y_k^{1/2}(d)\(1-\frac{k ^{1/2} \wt \vf(C|\De|)}{(U^*)^{1/2}}\)\)=1\nn. 
\eea
 
 We now show that this implies (\ref{rev.2quu}). 
 Let $  u,v,d\in \De$. 
We write, 
\begin{eqnarray}
\eta_{i}^{2}(u)- \eta_{i}^{2}(v)&=&( \eta_{i}(u)- \eta_{i}(v))( \eta_{i}(u)+ \eta_{i}(v))
\label{16.40u}\\
&=& ( \eta_{i}(u)- \eta_{i}(v))( 2 \eta_{i}(d)+( \eta_{i}(u)- \eta_{i}(d))+( \eta_{i}(v)- \eta_{i}(d)))   \nonumber.
\end{eqnarray} 
Consequently, 
\bea
&& \lim_{h\to 0}\sup_{\stackrel{|u-v|\le h }{ u,v\in\De}}\frac{Y_{k }(u)-Y_{k }(v) }{      \vf (u,v)   }\label{5.4mm2}\\&&\nn\qquad \ge \lim_{h\to 0}\sup_{\stackrel{|u-v|\le h }{ u,v\in\De}}\frac{2\sum_{i=1}^{k}( \eta _{i }(u)-\eta _{i }(v)) \eta_{i}(d)}{ \vf (u,v)} \\&&\qquad \qquad- \lim_{h\to 0}\sup_{\stackrel{|u-v|\le h }{ u,v\in\De}}\frac{ 2\sum_{i=1}^{k}( \eta _{i }(u)-\eta _{i }(v))( \eta _{i }(u)-\eta _{i }(d)) }{    \vf (u,v)}.\nn
  \eea
Using (\ref{rev.1wquqq})   we see that, 
 \bea 
 &&\lim_{h\to 0}\sup_{\stackrel{|u-v|\le h }{ u,v\in\De}}\frac{ 2\sum_{i=1}^{k}  | \eta _{i }(u)-\eta _{i }(v)||\eta _{i }(u)-\eta _{i }(d) | }{    \vf (u,v)} \\&&\qquad \le \lim_{h\to 0}\sup_{\stackrel{|u-v|\le h }{ u,v\in\De}}\frac{ 2\sum_{i=1}^{k}  | \eta _{i }(u)-\eta _{i }(v)|  } {    \vf (u,v)}\sup_{u\in \De}|\eta _{i }(u)-\eta _{i }(d) |\nn\\
&&\qquad  \le \nn 2k\sup_{1\le i\le k}\sup_{u\in \De}|\eta _{i }(u)-\eta _{i }(d) |:=\De^*_k. 
  \eea 
  Consequently,	
\bea
&& \lim_{h\to 0}\sup_{\stackrel{|u-v|\le h }{ u,v\in\De}}\frac{Y_{k }(u)-Y_{k }(v) }{      \vf (u,v)   }\label{5.4mm2}\\&&\nn\qquad \ge \lim_{h\to 0}\sup_{\stackrel{|u-v|\le h }{ u,v\in\De}}\frac{2\sum_{i=1}^{k}( \eta _{i }(u)-\eta _{i }(v)) \eta_{i}(d)}{ \vf (u,v)}  - \De^*_k.\nn
  \eea

   We now reverse the argument that leads to (\ref{2.50nn3}) to see that,\begin{eqnarray}
&& \lim_{h\to 0}\sup_{\stackrel{|u-v|\le h }{ u,v\in\De}} \frac{\sum_{i=1}^{k}( \eta _{i }(u)-\eta _{i }(v)) \eta_{i}(d)   }{    \   \vf (u,v)} \label{} \\
&&\qquad    \ge \lim_{h\to 0}\sup_{\stackrel{|u-v|\le h }{ u,v\in\De}} \frac{\sum_{i=1}^{k}\(V_{i }(v,d) \eta_{i}(u)-V_{i }(u,d)\eta_{i}(v)\) \eta_{i}(d)  }{     \vf (u,v)} \nn\\
&&\hspace{.75 in}- \lim_{h\to 0}\sup_{\stackrel{|u-v|\le h }{ u,v\in\De}}  \frac{| \sum_{i=1}^k \eta_i(d) G_{i}(u,v)| }{    \vf  (u,v)}, \qquad a.s.\label{}\nn
\end{eqnarray}
Consequently, using (\ref{1.51}) and  Lemma \ref{lem-3.1}, see (\ref{1.61}), we have   that for all $|\De|$ sufficiently small, 
\begin{eqnarray}
&& \lim_{h\to 0}\sup_{\stackrel{|u-v|\le h }{ u,v\in\De}} \frac{\sum_{i=1}^{k}( \eta _{i }(u)-\eta _{i }(v)) \eta_{i}(d)   }{    \   \vf (u,v)} \label{} \\
&&\qquad    \ge   Y^{1/2}_k(d)\(1-\frac{k^{^{1/2} }  \wt \vf(C|\De|)}{(U^*)^{1/2}}\) -\frac{k^{1/2} \wt \vf(C|\De|) }{(U^{*})^{1/2} }Y^{1/2}_k(d) \nn,\qquad a.s.,
\end{eqnarray}
which by (\ref{5.4mm2}) gives,
\be  \lim_{h\to 0}\sup_{\stackrel{|u-v|\le h }{ u,v\in\De}}\frac{Y_{k }(u)-Y_{k }(v) }{      \vf (u,v)   } \ge 2   Y_k^{1/2}(d)\(1-\frac{2k^{1/2}\wt \vf(C|\De|) }{(U^*)^{1/2}}\) - \De^*_k ,  \qquad a.s.
  \ee
 When $\De'\subset\De$,
\be  \lim_{h\to 0}\sup_{\stackrel{|u-v|\le h }{ u,v\in\De}}\frac{Y_{k }(u)-Y_{k }(v) }{      \vf (u,v)   } \ge  \lim_{h\to 0}\sup_{\stackrel{|u-v|\le h }{ u,v\in\De'}}\frac{Y_{k }(u)-Y_{k }(v) }{      \vf (u,v)   }. 
  \ee
Therefore  for any $\De'\subset\De$ with  $d\in\De'$,
\be  \lim_{h\to 0}\sup_{\stackrel{|u-v|\le h }{ u,v\in\De}}\frac{Y_{k }(u)-Y_{k }(v) }{      \vf (u,v)   } \ge 2   Y_k^{1/2}(d)\(1-\frac{2k^{1/2} \wt\vf(C|\De|')}{(U^*)^{1/2}}\) - (\De')_k^*,  \qquad a.s.\label{2.61}
  \ee
   for all $|\De'|$  sufficiently small.

   Let, 
\begin{equation} \label{}
\wt\si^{2}(x) =\sup_{1\le i\le k}\sup_{|u-v|\le x}E(\eta_{i}(u)-\eta_{i}(v))^2.
\end{equation}
 It follows from  \cite[Lemma 2.2]{MRECP}
that for all $|\De'|$ sufficiently small,
\begin{equation} \label{2.63a}
 P\((\De')^*\ge 4k^2\((1+2C')\wt\si^{2}(|\De'|)\log 1/|\De'|\)^{1/2}\)\le 4k^2|\De'|^{C'},
\end{equation} 
for any positive constant $C'$. 

  Let $\De'(d,n)\subset\De $ be an interval of length $1/n$ that contains $d$. Using  (\ref{2.61}) and (\ref{2.63a}) we see that,  
\bea  \lim_{h\to 0}\sup_{\stackrel{|u-v|\le h }{ u,v\in\De}}\frac{Y_{k }(u)-Y_{k }(v) }{      \vf (u,v)   } &\ge  &  2Y_k^{1/2}(d)\(1-\frac{ 2k^{1/2} \wt\vf(C/n)}{(U^*)^{1/2}}\) \\
&&  \nn  +\,4k^2\((1+2C')\wt\si^{2}(1/n)\log 1/n\)^{1/2},
   \eea
except possibly on a set of measure $4k^2/n ^{C'}$.   Taking $n\to \ff$  and using (\ref{2.38})   gives (\ref{rev.2qu6}) and consequently (\ref{rev.2qu5}) which is the lower bound in (\ref{rev.2quu}).

\medskip     
Exactly the same argument that   gives (\ref{2.61}) but with the inequalities reversed shows that for any $(m,n)$ and $\De_{n,m}=[\frac{m-1}{n},\frac{m+1}{n}]\in [0,1]$, for all $n$ sufficiently large,    and any $d_{m,n}\in \De_{n,m}$,
 \bea&&  \lim_{h\to 0}\sup_{\stackrel{|u-v|\le h }{ u,v\in [\frac{m-1}{n},\frac{m+1}{n}]}}\frac{Y_{k }(u)-Y_{k }(v) }{      \vf (u,v)   } \label{1.62}\\
 &&\qquad\le 2   Y_k^{1/2}(d_{m,n})\(1 +\frac{2k^{1/2}  \wt\vf( 2C/n)}{(U^*)^{1/2}}\) +(\De_{n,m} )_k^*, \qquad a.s.\nn 
  \eea

  Let $\wt \De_{n} =\{m\,|\, [\frac{m-1}{n},\frac{m+1}{n}]\cap \De \neq \emptyset\}$.  Then for any   $d_{m,n}\in \De_{n,m}\cap \De$ we have,
 \bea \label{2.63} &&\lim_{h\to 0}\sup_{\stackrel{|u-v|\le h }{ u,v\in\De}}\frac{Y_{k }(u)-Y_{k }(v) }{      \vf (u,v)   }\\
&&\qquad \leq \sup_{m\in \wt \De_{n} }\,\,\lim_{h\to 0}\sup_{\stackrel{|u-v|\le h }{  u,v\in [\frac{m-1}{n},\frac{m+1}{n}]}}\frac{Y_{k }(u)-Y_{k }(v) }{      \vf (u,v)   } \nn\\&&\qquad \le 2 \sup_{m\in \wt \De_{n} } \(  Y_k^{1/2}(d_{m,n})\(1 +\frac{2k^{1/2}  \wt \vf(  2C/n)}{(U^*)^{1/2}}\) + ( \De_{n,m})_k^*\) \nn\\ 
&&\qquad \le 2 \sup_{v\in   \De  } Y_k^{1/2}(v)\(1 +\frac{ 2k^{1/2}  \wt \vf(  2C/n)}{(U^*)^{1/2}}\) +2 \sup_{m\in \wt \De_{n} } ( \De_{n,m})_k^*  \nn. 
  \eea
  Using this and (\ref{2.63a}) we see that,
\bea \label{2.63} &&\hspace{-.15in}\lim_{h\to 0}\sup_{\stackrel{|u-v|\le h }{ u,v\in\De}}\frac{Y_{k }(u)-Y_{k }(v) }{      \vf (u,v)   }\\  &&\hspace{-.1in}  \le 2\sup_{v\in\De}  Y_k^{1/2}(v)\(1 +\frac{2k^{1/2} \wt \vf(  2C/n)}{(U^*)^{1/2}}\) +\,4k^2\((1+2C')\wt\si^{2}(2/n)\log 2/n\)^{1/2} \nn. 
  \eea
except possibly on a set of measure $4k^2n \(\frac{2}{n}\)^{C'}$. Taking $C'>2$ and letting $n\to \ff$ gives the upper bound in (\ref{rev.2quu}).

  \medskip  To prove (\ref{rev.2qusa}) we use essentially the same proof  that extends Theorem \ref{theo-1.3m} to Corollary \ref{cor-2.1}.  We start with   $\ov \eta_{i}=\{\eta_{i}(t)+\xi,t\in [0,1]\}$,   $i=0,\ldots,k$, where $\xi$ is a standard normal random variable independent of  $\eta_{i}=\{\eta_{i}(t),t\in [0,1]\}$,    $i=0,\ldots,k$. 
  But note that $E( \ov \eta_i^2(t))\geq 1$ so that we can apply (\ref{rev.2quu}) to obtain (\ref{rev.2qusa}).  Applying this with $s=0$ the shows that  (\ref{rev.2quu}) holds without requiring that
   $U^*>0$. \qed

\begin{lemma} \label{lem-3.1} Let $G_{i}(u,v)$ be as given in (\ref{1.45}).
Then  for all $|\De|$ sufficiently small,
 \begin{equation}
\lim_{h\to 0}\sup_{\stackrel{|u-v|\le h }{ u,v\in\De}}  \frac{| \sum_{i=1}^k a_i G_{i}(u,v)| }{    \vf  (u,v)}\le \frac{k \wt \vf( C|\De| ) }{(U^{*})^{1/2} }  ,\qquad a.s.\label{rev.ioi}
\end{equation} 
and,
\begin{equation}
\lim_{h\to 0}\sup_{\stackrel{|u-v|\le h }{ u,v\in\De}}  \frac{| \sum_{i=1}^k \eta_i(d) G_{i}(u,v)| }{    \vf  (u,v)}\le   \frac{ k\wt \vf(C|\De| ) }{(U^{*})^{1/2} }Y_{k}^{1/2}(d) ,\qquad a.s.\label{1.61}
\end{equation} 
\end{lemma}

\Proof   We first  show that,
\be  
 |G_{i}(u,v)|\leq   \frac{\si_{i}( d,v )}{U_{i}^{1/2}(d,d)} |\eta_{i} (u)- \eta_{i} (v)|+\frac{\si_{i}( u,v )}{U_{i}^{1/2}(d,d)} |\eta_{i} (v)|.\label{Gbound}
 \ee
We write,
\be \label{a02}
 G_{i} (u,v)=(1-V_{i}(v,d) )(\eta_{i} (u)-\eta_{i}(v)) +(V_{i}(u,d)-V_{i}(v,d))\eta_{i} (v) .
\ee
Note that,
\bea  
 |(V_{i}(u,d)-V_{i}(v,d))| &=&\frac{|\(E(\eta_{i} (u)- \eta_{i} (v)) \eta_{i} (d) \)|}{ U_{i}(d,d) }\\
 &\le&\frac{\(E(\eta_{i} (u)- \eta_{i} (v))^{2} E\eta_{i}^{2} (d)\)^{1/2}}{ U_{i}(d,d)}\nn\\
  &\le&\frac{\si_{i}( u,v )}{ U_{i}^{1/2}(d,d) }\nn,
\eea
and,
\bea 
 (1-V_{i}(v,d))  &=&\frac{\(E(\eta_{i} (d)- \eta_{i} (v)) \eta_{i} (d) \)}{ U_{i}(d,d) }\\
 &\le&\frac{\(E(\eta_{i} (d)- \eta_{i} (v))^{2} E\eta_{i}^{2} (d)\)^{1/2}}{ U_{i}(d,d)}\nn\\
  &\le&\frac{\si_{i}( d,v )}{ U_{i}^{1/2}(d,d)}  \nn.
\eea
These inequalities give (\ref{Gbound}). 

 It follows from (\ref{rev.1wquqq}) that, 
\be
\lim_{h\to 0}\sup_{\stackrel{|u-v|\le h }{ u,v\in\De}} \frac{ | \eta _{i }(u)-\eta _{i }(v)|   }{    \    \vf (u,v) } =1\label{rev3.4}.
\ee
  Using this  and (\ref{Gbound}) we have,  
\bea \label{2.60nn}
&&    \lim_{h\to 0}\sup_{\stackrel{|u-v|\le h }{ u,v\in\De}}\frac{   |G_{i}(u,v)|  }{   \vf (u,v)} \le    \lim_{h\to 0}\sup_{\stackrel{|u-v|\le h }{ u,v\in\De}}\frac{   |\eta_{i} (u)- \eta_{i} (v)|\si_{i}( d,v )  }{  U_{i}^{1/2}(d,d)\, \vf (u,v)}\,\,\,\\
 && \hspace{1.75in}+  \lim_{h\to 0}\sup_{\stackrel{|u-v|\le h }{ u,v\in\De}}\frac{    \si_{i}( u,v )   }{  U_{i}^{1/2}(d,d)\,  \vf (u,v)}\eta_{i}(v) \nn
 \eea
\bea 
 &&\hspace{.5in} \nn\le \sup_{d, v\in \De}\frac{\si_i( d,v )}{U_i^{1/2}(d,d)}+ \lim_{h\to 0}\sup_{\stackrel{|u-v|\le h }{ u,v\in\De}}\frac{   \si_{i}( u,v )   }{   \vf (u,v)}\frac{\sup_{ v\in \De}|\eta\(v\)|}{U_i^{1/2}(d,d)} \eea

  It follows from (\ref{rev3.4}) and  \cite[Lemma 7.1.1]{book} that,
\begin{equation} \label{2.72}
  \lim_{h\to 0}\sup_{\stackrel{|u-v|\le h }{ u,v\in\De}}\frac{    \si_{i}( u,v )   }{   \vf (u,v)} =0. 
\end{equation}
Using this in the last line of (\ref{2.60nn}) we see that for all $|\De|$ sufficiently,
\begin{equation} \label{2.73}
  \lim_{h\to 0}\sup_{\stackrel{|u-v|\le h }{ u,v\in\De}}\frac{   |G_{i}(u,v)|  }{   \vf (u,v)} \le \sup_{d, v\in \De}\frac{\si_i( d,v )}{U_i^{1/2}(d,d)}\le  \frac{\wt \vf(C|\De| )}{(U^ *)^{1/2}},  \qquad a.s. 
\end{equation}

Note that since all $|a_i|$  are less than or equal to 1 and     $\eta_i(d)\leq \|(\eta_1(d),\ldots, \eta_k(d))\|_2=Y^{1/2}_k(d) $ for each $i$,
\begin{equation} \label{}
 \Big| \sum_{i=1}^k a_i G_{i}(u,v)\Big| \le   \sum_{i=1}^k  |G_{i}(u,v)| 
\end{equation}
and,
\begin{equation} \label{}
 \Big| \sum_{i=1}^k \eta_i(d) G_{i}(u,v)\Big| \le  Y_k^{1/2}(d) \sum_{i=1}^k  |G_{i}(u,v)| .
\end{equation}
 Using  these two inequalities and (\ref{2.73}) we get (\ref{rev.ioi})
 and (\ref{1.61}).

\begin{remark} {\rm \label{rem-3.2} In this paper we apply Theorem \ref{theo-2.1m} when $G_{r,s}$ is as given in (\ref{2.21p}). That is $\{\eta_i \}_{i=1}^{r-1}$ are independent copies   $\{\eta_{0, i } \}_{i=1}^{r-1}$ of $\eta_0$, so that,  \begin{equation}
\sum_{i=1}^{r-1} \eta _{0, i }\,\,a_{i}  \stl \eta_{0}  \(\sum_{i=1}^{r-1}a^{2}_{i}\)^{1/2}.\label{}
\end{equation}
Therefore, the condition in (\ref{rev.1wquqq}) that 
   defines $ \GG_{\vf,\De,k}$ can be simplified and written as,
   \be
\lim_{h\to 0}\sup_{\stackrel{|u-v|\le h }{ u,v\in\De}} \frac{ ( \eta _{0 }(u)-\eta _{0 }(v)) a_{0}  + ( \eta _{r }(u)-\eta _{r }(v)) a_{r} }{    \    \vf (u,v) } =1 \label{} ,
\ee
  for all   intervals $\De$ in $[0,1]$ and all     $a_0^2+a_r^2=1 $.  

}\end{remark}

\subsection{Local and uniform moduli of continuity of Gaussian processes with stationary increments}

The next theorem pulls together several results that are scattered in  \cite[Sections 7.2 and 7.3]{book}.

  \bt\label{theo-2.3} Let $G=\{G(x),x\in [-1,1]\}$ be a Gaussian process
with stationary increments with the property that,
\begin{equation} \label{5.13nn}
\si   ^2(x-y)=E(G(x)-G(x))^2 =\frac{1}{\pi}\int_{-\ff}^\ff\frac{1- \cos\la (x-y)}{  \th(\la)}\,d\la,\qquad x,y\in R^{1},  
\end{equation}
where   $\th(\la)$ is positive and symmetric and 
\begin{equation} \label{3.82}
 \int_{-\ff}^{\ff} \frac{1\wedge\la^2}{\th(\la)}\,d\la<\ff.
\end{equation}
Then when $\th$ is a regularly varying function at infinity with index $1<\ga\le 2$,
 \be
\lim_{\de\to 0}\sup_{ {|u| \le\de}  }
\frac{|G(u+d)-G(d)|}{ (2\si^2(u)\log\log 1/ |u|)^{1/2}}= 1\qquad \mbox {a.s.}
\label{5.36}
\ee
  for all $d\in [-1,1]$ and,
\be
\lim_{\de\to 0}\sup_{\stackrel{|u-v|\le\de}{u,v\in [-1,1]} }
\frac{|G(u)-G(v)|}{ (2\si^2(|u-v|)\log 1/|u-v|)^{1/2}}= 1\qquad \mbox {a.s.} 
\label{5.37p}
\ee

    \et

\Proof    The lower bound in (\ref{5.36}) follows from  \cite[Theorem 7.3.2]{book}. The upper bound in (\ref{5.36}) follows from  \cite[Corollary 7.2.3]{book} using the observations in    \cite[ Remark 7.3.4]{book}.  

The statement in (\ref{5.37p}) follows from  \cite[Theorem 7.3.3 and Remark 7.3.4]{book}.
  Note that the condition in the second paragraph of  \cite[Theorem 7.3.3]{book} is satisfied when (\ref{5.13nn}) and (\ref{3.82}) are satisfied. 

(In  \cite[(7.223)]{book} include the condition that $u,v\in [0,1] $ as in \cite[(7.228)]{book}.)
 \qed

\begin{remark} {\rm \label{} The results in  
 Theorem \ref{theo-2.3} also hold when $\th$ is
  a regularly varying function at infinity with index 1 and some additional regularity conditions hold. We leave this case to the interested reader.
  
  Also note that the title of  \cite[Section 7.3]{book} is misleading. In this section we only require that that the increments variance of the Gaussian process with stationary increments can be written as in 
\cite[(7.190)]{book}. 
  
  }\end{remark}

\section{Local and Uniform Moduli of Continuity for the local times of Rebirthed Markov Processes.}\label{sec-app} 

We use Theorem \ref{theo-ITcond} to obtain local and uniform moduli of continuity for the local times of rebirthed Markov processes.   First we greatly simplify the conditions in (\ref{2.21p}) and (\ref{70.6n}) and then we use   generalizations of the techniques in   \cite[Theorem 9.5.1]{book} to complete the proof.

It seems useful to summarize the definitions in the beginning of this paper. Let   $\YY=(\Om,  \FF_{t},  \YY_t,\th_{t},   P^x)$ be a transient symmetric Borel right process with state space $S$  and  continuous strictly positive  $p-$potential   densities  
\be
  u^{p}=\{ u^{p}(x,y), x,y\in S \},\quad p\geq 0,
\ee
 with respect to some $\si$--finite  positive measure $m$ on $S$. Let $\wt Z$ be a  rebirthed Markov process  of $\YY$ as described in the beginning of this paper and let  $\wt L_{t}^{y}$ denote the local time of $\wt Z$
  normalized so that,
\begin{equation}
  \wt E^{ x}\(\int_{0}^{\ff}  e^{-ps}\, d_{s}\wt L^{y}_s\)=   w^p(x,y),\qquad x,y\in S, \label{int.3q}
\end{equation}
where $w^p(x,y)$ is given in (\ref{int.1}).

\medskip  We assume that $S\subseteq R^1$. In the examples we give in Section \ref{sec-examp} the state space $S$ is either $R^1$, $R_+$, $R^1-\{0 \}$ or $R_+-\{0  \}$.

  \bt\label{theo-localmod} 
 Let $\eta_{0}=\{\eta_{0} (x), x\in [0,1]\cap S\}$ be a Gaussian process with covariance $ u^{0}(x,y)$ and  $\eta_{p}=\{\eta_{p} (x), x\in [0,1]\cap S\}$ be a Gaussian process with covariance $ u^{p}(x,y)$. 
For $p\ge 0$ set, 
  \begin{equation} \label{4.4}
  (\si_{p})^{2} (x,y)= u^{p}(x,x)+ u^{p}(y,y)-2  u^{p}(x,y). 
\end{equation} 
 If for $d\in [0,1]\cap S$,    
   \begin{equation}
 \lim_{x,y\to d } \frac{\si_{p }^2(x,y)}{\si_{0} ^2(x,y)}=1\label{1.4amm},
   \end{equation}
  and   $\eta_{0} 
$ has a local modulus of continuity at   $d $   of the form,
  \begin{equation} \label{1.7m}
   \limsup_{x\to 0} \frac{|\eta_{0} (d+x)-\eta_{0}  (d)|} {   \phi(x)} =  1    \qquad a.s.
\end{equation}
for some increasing function $\phi$, then 
   \be
   \limsup_{x\to 0}\frac{|\wt L_{t}^{x+d}-\wt L_{t}^{d}|}{\phi(x)  }=  2\(\wt L_{t}^{d}\)^{1/2},  \quad a.e.\,\,\, t,  \,\,\, \wt P^{y}\,\,
a.s.\label{rev.2}
\ee
for all $y\in S$.

\et

\bt\label{theo-unifmod}   Let $\eta_{0} (x)$ and $\eta_{p} (x)$ be    as described in Theorem \ref{theo-localmod} and assume that (\ref{2.38}) and (\ref{3.39}) hold.       When,
 \be
 \lim_{h\to 0}\sup_{\stackrel{|u-v|\le h }{ u,v\in\De}} \frac{ ( \eta _{0 }(u)-\eta _{0}(v)) a_{0} + (\eta _{p }(u)-\eta _{p}(v)) a_{p} }{    \    \vf (u,v) } =1 \label{rev.1wqi} ,
\ee
 for all   intervals $\De$ in $[0,1]\cap S$ and all     $\{a_0,a_p\}$  with  $a_0^2+a_p^2=1$     then,   \be
\lim_{h\to 0}\sup_{\stackrel{|u-v|\le h }{ u,v\in\De}} \frac{|\wt L_{t}^{u}-\wt L_{t}^{v}|}{\vf (u,v)  }= \sup_{u\in \De}  \(2\wt L_{t}^{u}\)^{1/2},  \quad a.e.\,\,\, t ,  \,\,\, \wt P^{y}\,\,
a.s.,\label{urev.2}
\ee
for all $y\in S$.  
\et 

 \noindent\textbf{Proof of Theorem \ref{theo-localmod} }   
 Let $F( \mathcal{C})$  be a set of  functions on a countable dense   subset $\mathcal{C}$ of $[0,1]\subseteq S$ and set,  
\be
B =\Big\{g\in F( \mathcal{C})\,\Big |  \,\limsup_{x\in \mathcal{C},   x\to 0}\frac{|g (d+x)-  g(d)|}{\phi(x)}=  \(2 {g(d)}\)^{1/2} \Big\}.\label{A.local}
\ee
  Now consider,
 \be G_{r,s}\(x\)=\sum_{i=1}^{r-1} { 1\over 2}(\eta_{i, 0 }(x)+s)^2+{ 1\over 2} (\eta_{p}(x)+s)^2. \label{2.21}
\ee
where $\{ \eta_{i, 0 },1\le i\le r-1\}$ are independent copies of $\eta_{ 0 }$. It follows from  (\ref{1.4amm}) and (\ref{1.7m}) and  Corollary \ref{cor-2.1} that,
\begin{equation} 
 \limsup_{t\to 0} \frac{ G_{r, s}\( d+t \)- G_{r, s}\(d\)} {   \phi(t)} =   2 \(  G_{r, s}\(d\)\)^{1/2} , \label{1.35ksa}
\ee  almost surely. Consequently,
 for each $r\geq 1$,   
\be
  P_{G_{r, s}} \(G_{r, s} \in B\)=1.\label{az.2p}
\ee 
It now follows from   Theorem \ref{theo-ITcond},  conditional on $\ze_{r-1}<\la \leq \ze_{r}$, that for almost all $\om'\in
\Om_{G_{r,s}}$, the probability space of $G_{r,s}$, 
\be  \limsup_{ x\in \mathcal{C},   x\to 0}
\Bigg|\frac{\wt  L^{d+x}_{\la}-\wt  L^{ d}_{\la}}{\phi(x)}\!+  \frac{   G_{r,s}( d+x,\om')- G_{r,s} (d,\om') }{\phi(x)} \Bigg| =
  \,  \({2( \wt  L^{ d}_{\la}+ G_{r,s} (d,\om')}\)^{1/2}, \label{az.7}  
\ee
 $ \wt P_{\la}^{y}$ almost surely, for all $y\in S$. 

  Using (\ref{az.2p}) and (\ref{az.7}) we get that conditional on $\ze_{r-1}<\la \leq \ze_{r}$,
\be 
  \bigg|  \limsup_{x\to 0}\frac{\wt L^{ d+x}_{\la}-\wt L^{d}_{\la}}{\phi(x)  }-\(2\wt L^{d}_{\la}\)^{1/2}\bigg|\le 
  2\(       G_{r,s} (d,\om')\)^{1/2}, \quad   \wt P_{\la}^{y}\,\,\, a.s.  \nn
\ee
 By (\ref{2.21}),
 \be G_{r,s}\(d,\om'\)=\sum_{i=1}^{r-1} { 1\over 2}(\eta_{i, 0 }(d,\om')+s)^2+{ 1\over 2} (\eta_{p}(d,\om')+s)^2. \label{2.21q}
\ee
Since the random variables $|\eta_{0,i }(d)|$, $i=1,\ldots,r-1$ and $|\eta_{p}(d)|$ take values arbitrarily close to 0 with probability greater than 0, and we can take $s$ arbitrarily close to 0, we see that 
\begin{equation}
\wt P_{\la}^{y}\(\limsup_{ x\in \mathcal{C},   x\to 0}
\Bigg|\frac{\wt   L^{d+x}_{\la}- \wt  L^{ d}_{\la}}{\phi(x)} \Bigg|=      \(  2\wt   L^{ d}_{\la}\)^{1/2}\,\Big |\,  \ze_{r-1}<\la \leq \ze_{r} \)=1.\label{az.9}
\end{equation}
Since this is true for all $r\geq 1$ it follows from the law of total probability that, for all $y\in S$,
\begin{equation}
\wt P_{\la}^{y}\(\limsup_{ x\in \mathcal{C},   x\to 0}
\Bigg|\frac{\wt   L^{d+x}_{\la}- \wt  L^{ d}_{\la}}{\phi(x)} \Bigg|=      \( 2 \wt   L^{ d}_{\la}\)^{1/2}\)=1.\label{az.9a}
\end{equation}
Equivalently, by Fubini's theorem, 
\begin{equation}
\int_{0}^{\ff} \wt P^{y}\(\limsup_{ x\in \mathcal{C},   x\to 0}\Bigg|\frac{\wt   L^{d+x}_{t}- \wt  L^{ d}_t}{\phi(x)} \Bigg|=      \( 2 \wt   L^{ d}_{t}\)^{1/2}\)\, pe^{-pt}\,dt=1,\label{az.9a2}
\end{equation}
 for all $y\in S$. Since,
\begin{equation}
\int_{0}^{\ff}  pe^{-pt}\,dt=1,\label{2.28}
\end{equation} we have, 
\begin{equation}
  \wt P^{y}\(\limsup_{ x\in \mathcal{C},   x\to 0}\Bigg|\frac{\wt   L^{d+x}_{t}- \wt  L^{ d}_t}{\phi(x)} \Bigg|=      \( 2\wt   L^{ d}_{t}\)^{1/2}\)=1,\label{az.9a3}
\end{equation}
for almost all $t>0$, for all $y\in S$.

\medskip It follows from by Lemma \ref{lem-ltdecomp} that  $\wt  L^{x}_t$ is continuous    in $x$. Therefore, we can remove the condition that $x\in \mathcal{C}$. This gives us  (\ref{rev.2}).
\qed

 \noindent\textbf{Proof of Theorem \ref{theo-unifmod} }  The proof is essentially  the same as the proof Theorem \ref{theo-localmod} except  we now use Theorem \ref{theo-2.1m} and Remark \ref{rem-3.2} 
and  instead of $B$ in (\ref{A.local}) we  use,
\be
A =\Big\{g\in F( \mathcal{C})\,\Big |  \,\lim_{h\to 0}\sup_{\stackrel{|u-v|\le h }{ u,v\in\De\cap \mathcal{C}}}\frac{|g (u)-  g(v)|}{\vf(u,v)}= \sup_{u\in \De} \(2{g(u)}\)^{1/2} \Big\},\label{B.unif}
\ee
and  the fact that the random variables $\sup_{u\in \De}|\eta_{i }(u)|$, $i=1,\ldots,r-1$ and $\sup_{u\in \De}|\eta_{p}(u)|$ take values arbitrarily close to 0 with probability greater than 0.    
\qed
 
   \section{Moduli of continuity of the local of  times of certain fully rebirthed Markov processes}\label{sec-examp}
 
 \begin{example} \label{}{\rm    \textbf{L\'evy Processes } 
 Let $Z=\{Z_{t};t\ge 0 \}$ be a real valued symmetric 
L\'evy process with characteristic exponent $\psi$, i.e., 
\begin{equation} \label{1.27nn}
  E\(e^{i\la Z_{t}}\)=e^{-t\psi(\la)},\qquad \forall\,t\ge 0,
\end{equation}
where,
\begin{equation} \label{1.20nn}
  \frac{1}{ 1+\psi(\la)}\in L^{1}(R^1).
\end{equation}
Note that the condition that   $Z$ is symmetric implies that $\psi$ is an even function.

 For $\bb>0$, the $\bb-$potential density of $Z$ with respect to Lebesgue measure  is,  
\begin{equation} \label{1.21nn}
u^{\bb}:=u^{\bb}(x,y)=u^{\bb}(x-y)=\frac{1}{2\pi}\int_{-\ff}^{\ff}\frac{ \cos\la (x-y)}{\bb+ \psi(\la)}\,d\la,\qquad x,y\in R^{1};  
\end{equation}
see e.g.,   \cite[(4.84)]{book}.

 We define
 \be\label{5.4}
  {\frak u}^{0}(x,y)=\phi(x)+\phi(y)-\phi(x-y),
\end{equation}
where \begin{equation} \label{lv.22}
\phi(x)=\frac{1}{2\pi}\int_{-\ff}^{\ff}\frac{1-\cos\la x}{ \psi(\la)}\,d\la ,\qquad x \in R^{1}.  
\end{equation}
 (Note that   $ \frak {u}^{0}(x,y)$ is not $ {u}^{\bb}(x,y)$ with $\bb=0$.) 

For $\bb> 0$ set,
  \bea \label{4.4pq}
  (\si^{\bb})^{2} (x-y)&=& u^{ \bb}(x,x)+ u^{ \bb}(y,y)-2  u^{\bb}(x,y)\\
   &=& \nn\frac{1}{\pi}\int_{-\ff}^\ff\frac{1- \cos\la (x-y)}{\bb+ \psi(\la)}\,d\la,
\eea 
and 
 \bea \label{4.4pqq}
  (\si^{0})^{2} (x-y)&=& \frak u^{0}(x,x)+ \frak u^{ 0}(y,y)-2  \frak u^{0}(x,y)\\
   &=& \nn\frac{1}{\pi}\int_{-\ff}^\ff\frac{1- \cos\la (x-y)}{  \psi(\la)}\,d\la.
\eea
Note that,  
\begin{equation}
( \si ^0)^2(x-y) =2\phi(x-y) .\label{uisv.2}
\end{equation} 

We show in  \cite[Lemma 8.1]{MRLIL} that,
 for all $\bb,\bb'\ge0$,
 \begin{equation}
  \lim_{x \to 0} \frac{(\si ^\bb)^2(x )}{ (\si ^{\bb'})^2(x )}=1\label{5.21}.
   \end{equation}
   In addition it follows from  \cite[Lemma 7.3.1]{book} that when $\psi$ is regularly varying at infinity with index $1<r\le 2$,   for all $\bb\ge 0$ 
\begin{equation} \label{5.22}
  (\si^\bb)^2(x)\sim C_r\frac{1}{x\psi(1/x)},\qquad \text{as $\,x\to 0$},
\end{equation}
where,
\begin{equation} \label{}
  C_r=\frac{4}{\pi}\int_{0}^\ff\frac{\sin^2 s/2}{s^r}\,ds.
\end{equation}
 
Finally, we note that by  \cite[Section  II Theorem 19]{Bertoin}, for all $\bb> 0$,
\begin{equation} \label{5.11}
  u^\bb(x)>0,\qquad \forall\, x\in R^1.
\end{equation}

\medskip \noindent \textbf{Case 1.} \label{37}  Let   $Y_1=\{Y_{1,t},t\ge 0 \}$ be the process with state space $S_{1}=R^{1}$ that is obtained by killing  $Z$ at the end of an independent exponential time with mean $1/\bb$. In this case    $Y_1$ has $0-$potential density $u^{\bb}$ and 
   $p-$potential density   $u^{\bb+p}$.  Since  $P^{y}\(Z_{t}\in R^{1}, \forall t\geq 0\)=1$ for all $y\in R^{1}$, we have,\begin{equation} \label{5.4mm}
  \int_{-\ff} ^\ff u^{p}(x,y)\,dy=\frac1{p}, \qquad \forall\, p>0.
\end{equation} 
    It follows from (\ref{193.3a}) that
 in this case the function $w^p(x,y)$ in (\ref{int.1}) is
\begin{equation} \label{int.2simple}
  w_1^{p}(x,y)=u^{\bb+p}(x-y)+ \frac{\bb}{p}f_1(y),
\end{equation}
where \begin{equation} \label{int.2}
  f_{1}(y)=\int_{ -\ff} ^{\ff}  u^{\bb+p}(x-y)\,d\mu(x).
\end{equation}
For more details see  \cite[(7.29)]{MRLIL}. (This may be confusing but here the 0--potential is $u^\bb$.)

\medskip\noindent \textbf{Case 2.}  We consider the process  $Y_{2}=\{Y_{2,t},t\ge 0 \}$ with state space $S_{2}=R^{1}-\{0\}$ that is obtained by killing $Y_1$ when it hits 0. The $0-$potential density of $Y_2$ is,
\begin{equation} \label{5.6}
 v^{\bb} :=v^{\bb}(x,y)=u^\bb(x-y)-\frac{u^\bb(x )u^\bb( y)}{u^\bb(0)};
\end{equation} 
see  \cite[4.165]{book}.

The  $p-$potential density of     $Y_2 $   is  $v^{\bb+p}$.    Using (\ref{193.3a}) and (\ref{5.4mm}) we see that  the function $w^p(x,y)$ in (\ref{int.1}) is,
\begin{equation} \label{5.7}
  w_2^{p}(x,y)=v^{\bb+p}(x,y)+\frac{1}{\bb+p} \(\frac{\bb}{p}+\frac{u^{\bb+p}(x)}{u^{\bb+p}(0)}\)\frac{f_2(y)}{\|f_2\|_1},
\end{equation}
where $\|\cdot\|_1$ is the $L^1$ norm and
\begin{equation} \label{5.8}
  f_{2}(y)=\int_{ -\ff} ^{\ff}  v^{\bb+p}(x,y)\,d\mu(x).
\end{equation}
 For more details see  \cite[(7.33)]{MRLIL}.

\medskip\noindent \textbf{Case 3.}  Let $\bar Z=\{\bar Z_{t};t\ge 0 \}$ be a real valued symmetric 
L\'evy process with characteristic exponent $\psi$ 
where, 
\begin{equation} \label{1.20nn}
  \frac{1\wedge \la^2}{ \psi(\la)}\in L^{1}(R^1).
\end{equation} 
   and assume that 0 is recurrent for $\bar Z$,  or equivalently that   $u^{\bb}(0) $ in (\ref{1.21nn}) is infinite when $\bb=0.$ 
 Let  $Y_{3}=\{Y_{3,t},t\ge 0 \}$ be  a process with state space $S_{3}=R^{1}-\{0\}$   obtained by killing $\bar Z$ when it hits 0.    
   We show in \cite[Theorem 4.2.4]{book} that $ \{ \frak{u}^{0}(x,y),x,y\in R^1-\{0 \} \}$ is  the $0-$potential density of $Y_3$.  
   We show in   \cite[Example 4.5.4]{book} that for any $p>0$ the $p-$potential density   of     $Y_3 $ is $v^p(x,y)$ given in (\ref{5.6}).

  For these processes the function $w^p(x,y)$ in (\ref{int.1}) is given by (\ref{5.7}) and (\ref{5.8}) with $\bb=0$, i.e.  
 \begin{equation} \label{5.7mm}
  w_3^{p}(x,y)=v^{ p}(x,y)+\frac{1}{p} \,\,\frac{u^{p}(x)}{u^{p}(0)} \frac{f_3(y)}{\|f_3\|_1}.
\end{equation}
where \begin{equation} \label{5.13}
  f_{3}(y)=\int_{ -\ff} ^{\ff}  v^{ p}(x,y)\,d\mu(x).
\end{equation}
(This is (\ref{5.7}) and (\ref{5.8}) with $\bb=0$.)
 }\end{example}

\medskip Let $\{\wt Y_{i},i=1,2,3 \}$ denote the 
rebirthing of the process $\{  Y_{i},i=1,2,3 \}$ with probability measures $\mu$ supported on $S_{i}$. Let 
 $\wt L_{i,t}^{y}$ denote the local time of $\wt Y_i$
  normalized so that,
\begin{equation}
  \wt E^{ x}\(\int_{0}^{\ff}  e^{-ps}\, d_{s}\wt L^{y}_{i,s}\)=   w_i^p(x,y),\qquad x,y\in S_{i}. \label{5.12}
\end{equation}

Note that the examples in Cases
1., 2., and 3. are completely  determined by two quantities. The probability measure $\mu$ and the L\'evy exponent $\psi$. In our applications of Theorems \ref{theo-localmod} and \ref{theo-unifmod} we require that $\psi$ is regularly varying at infinity. It follows from  \cite[Lemma 7.4.10]{book} that for every  function $g$ that regularly varying at infinity with index $1<p<2$, there exists a L\'evy exponent $\psi$ such that,
\begin{equation} \label{}
  \lim_{\la\to \ff}\frac{\psi(\la)}{g(\la)}=1.
\end{equation}
In particular, $\psi$ is regularly varying at infinity with index $1<p<2$.

 \begin{theorem} \label{theo-5.1} Consider the examples in  Cases
1., 2., and 3.   determined by a  probability measure $\mu$ and   L\'evy exponent $\psi$. Assume that $\psi$ is regularly varying at infinity with index $1<p\le 2$. Then for $i=1,2,3$  and   $d\in S_{i}$,
   \be
   \limsup_{x\to d}\frac{|\wt L_{i,t}^{x}-\wt L_{i,t}^{d}|}{(2(\si^0)^2( x-d)\log\log  1/| x-d|)^{1/2}}= \(2\wt L_{i,t}^{d}\)^{1/2},  \quad a.e.\,\,\, t,  \,\,\, \wt P^{y}\,\,
a.s. \label{rev.2n}
\ee
for all $y\in S_{i}$. Furthermore,
  for any interval $\De$ in $S_{i}\cap [0,1]$,
   \be
\lim_{h\to 0}\sup_{\stackrel{|u-v|\le h }{ u,v\in\De}} \frac{|\wt L_{i,t}^{u}-\wt L_{i,t}^{v}|}{  (2(\si^0)^2(u-v) \log  1/| u-v|)^{1/2} }= \sup_{u\in \De}\(2\wt L_{i,t}^{u}\)^{1/2},  \quad a.e.\,\,\, t,  \,\,\, \wt P^{y}\,\,
a.s.\label{urev.2}
\ee
for all $y\in S_{i}$.  
\end{theorem}

\medskip \noindent\Proof   
  We first consider the local modulus of continuity. 

\medskip\noindent\textbf{Case 1.} 
    Let $G=\{G(x),x\in R^1 \}$ be a Gaussian process with covariance $u^{\bb}(x,y)$, $\bb>0$, as in (\ref{1.21nn}). It follows from Theorem \ref{theo-2.3}  and the hypothesis that $\psi$ is regularly varying at infinity with index $1<\ga\le 2$, that for all   $d\in R^1$,
\be
\lim_{\de\to 0}\sup_{ {|u| \le\de}  }
\frac{|G(u+d)- G(d)|}{ (2(\si^\bb)^2(u)\log\log 1/ |u|)^{1/2}}= 1\qquad \mbox {a.s.}
\label{5.24}
\ee  
We see  by (\ref{5.22}) that $(2\si^\bb)^2(u)\log\log 1/|u|$ is a regularly varying function at 0 with index $0<p'\le  1$. Therefore it is asymptotic at 0 to an increasing function, say $(\phi')^2$, so that,  
\be
\lim_{\de\to 0}\sup_{ {|u| \le\de}  }
\frac{|G(u+d)-G(d)|}{   \phi' (u) }= 1\qquad \mbox {a.s.}.
\ee
Using 
 Theorem \ref{theo-localmod} with $\eta_{0}=G$, a Gaussian process with covariance $u^{\bb}(x,y)$, and $\eta_{p}$ a Gaussian process with covariance $u^{\bb+p}(x,y)$,   and using (\ref{5.21}) we get   (\ref{rev.2n}) with   $ ((2\si^0)^2(u)  \log\log 1/|u|)^{1/2}$ replaced by 
 $(\phi')(u)$. However since these two functions are asymptotic at 0 we get (\ref{rev.2n}) as stated.

\begin{remark} \label{rem-5.1}{\rm \label{} In most of the examples in this paper we use functions asymptotic to  $(\si^0)^2$ at zero and require that they are regularly varying at zero with positive index. Consequently  they are asymptotic to an increasing function at 0. Using the argument just above   we see that we only need that the functions $\phi$ and $\vf$ in (\ref{1.7}) and (\ref{rev.1wquqq}) are asymptotic to an increasing function at zero.

}\end{remark}

  \medskip\noindent\textbf{Case 2.}
For $p>0$ set
\begin{equation} \label{4.4qq}
  (\si^p_{v})^{2} (x,y)=  v^p(x,x)+ v^p(y,y)-2   v^p(x,y), 
\end{equation} 
so that, 
\bea \label{5.14}
  (\si^p_{v})^{2} (x+d,y+d)&= & (\si^p )^{2} (x-y)-\frac{(u^p(x+d)-u^p(y+d))^2}{ 
   u^p(0)} \\&= & (\si^p )^{2} (x-y)-\frac{( (\si^p )^{2}(x+d)- (\si^p )^{2}(y+d))^2}{4 u^p(0)}.\nn
  \eea
We show in  \cite[(7.270)]{book}  
that for $p>0$,
\begin{equation} \label{5.15}
  |(\si^p )^{2}(x+d)- (\si^p )^{2}(y+d)|\le \frac{2u^p(0)}{u^p(x+d)+u^p(y+d)}(\si^p )^{2} (x-y),\qquad x,y\in R^{1}.
\end{equation}
Considering this and (\ref{5.14}) we see that for all $p>0$,   
 \begin{equation}
 \lim_{x,y\to 0} \frac{(\si_{v }^p)^2(x+d,y+d)}{ (\si ^p)^2(x-y)}=1\label{5.20v}.
 \ee
Therefore, it follows from (\ref{5.21}) that for all $p>0$, 
\begin{equation}
 \lim_{x,y\to 0} \frac{(\si_{v }^p)^2(x+d,y+d)}{ (\si ^0)^2(x-y)}=1\label{5.20vu}.
 \ee
 
  Using 
 Theorem \ref{theo-localmod} with $\eta_{0}=G$, a Gaussian process with covariance $u^{\bb}(x,y)$, and $\eta_{p}$ a Gaussian process with covariance $v^{\bb+p}(x,y)$,   and using   (\ref{5.20vu}) and Remark \ref{rem-5.1} we get (\ref{rev.2n}) as in the proof of Case 1.  
 
 \medskip\noindent\textbf{Case 3.}  We use  Theorem \ref{theo-localmod} with $\eta_{0}=G$, a Gaussian process with covariance $\frak{u}^{0}(x,y)$, and $\eta_{p}$ a Gaussian process with covariance $v^{p}(x,y)$. Then this follows as in Case 2 since we have (\ref{5.20vu}).  
 
  \medskip We now consider the uniform modulus of   continuity.

\medskip\noindent\textbf{Case 1.}    Let   $\eta_{\bb}=\{\eta_{\bb}(x),x\in R^1 \}$ be a Gaussian process with covariance $u^{\bb}(x,y)$ and $\eta_{\bb+p}=\{\eta_{\bb+p}(x),x\in R^1 \}$ be an independent Gaussian process with covariance $u^{\bb+p}(x,y)$ and as in   (\ref{1.21nn}).  Let,
\begin{equation} \label{}
  G(x)=a_0\eta_{\bb}(x)+a_1  \eta_{\bb+p}(x),\qquad \text{where}\,\,a_0^2+a_1^2=1.
\end{equation} 
Then
\begin{equation} \label{}
\si^2_G(x-y):= E(G(x)-G(y))^2= a_0^2(\si^\bb)^2(x-y)+a_1^2(\si^{\bb+p})^2(x-y)
\end{equation}
It follows from (\ref{5.22}) that $\si^2_G(x)$
is a regularly varying function at 0 with index $r-1$  when $\psi(x)$ is a regularly varying function at infinity   with index $r$. Consequently, it follows from Theorem \ref{theo-unifmod}  that,
\be
\lim_{\de\to 0}\sup_{\stackrel{|u-v|\le\de}{u,v\in [-1,1]} }
\frac{|G(u)-G(v)|}{ (2\si_0^2(|u-v|)\log 1/|u-v|)^{1/2}}= 1\qquad \mbox {a.s.} 
\label{5.37}
\ee
Using this in Theorem \ref{theo-unifmod} we get (\ref{urev.2}) in this case.

   \medskip\noindent\textbf{Case 2.}    Let   $\eta_{p}=\{\eta_{p}(x),x\in R^1 \}$ be a Gaussian process with covariance $u^{p}(x,y)$ and $\ov \eta_{ p}=\{\ov \eta_{ p}(x),x\in R^1 \}$ be a Gaussian process with covariance $v^{ p}(x,y)$.     
   We see by (\ref{5.6}) that,
 \begin{equation} \label{5.26}
 \ov \eta_{p}(x)=\eta_{p}(x)-\frac{u^p(x)}{u^p(0)}\eta_{p}(0). 
\end{equation}
Consequently,
\begin{equation} \label{5.26a}
  \eta_p(x)= \ov  \eta_{p}(x)+\frac{u^\bb(x)}{u^\bb(0)}\eta_\bb(0)\stl  \ov  \eta_{p}(x) +\frac{u^\bb(x)}{(u^\bb(0))^{1/2}}\xi,
\end{equation}
where $\xi$ is a standard normal random variable  that is independent of $\{\ov \eta_{p}(x),x\in R^1 \}$ and,
\begin{equation} \label{5.26q}
     \eta_{p} (u)-      \eta_{p} (v)\stl   \ov  \eta_{p}(u) -  \ov  \eta_{p}(v) +\frac{u^p(u)-u^p(v)}{({u^p(0)})^{1/2}}\xi ,
\end{equation}
Using  (\ref{5.15}) and (\ref{5.11}) we have, 
\bea \label{5.28}
 \bigg |\frac{u^p(u)-u^p(v)}{(u^\bb(0))^{1/2}}\xi \bigg |&\le& 2 \bigg |\frac{(\si^p )^{2}(x)- (\si^p )^{2}(y)}{(u^p(0))^{1/2}}\xi \bigg |\nn\\\
 &\le& \frac{4 (u^p(0))^{1/2}}{u^p(u)+u^p(v)}(\si^p )^{2} (u-v)|\xi  |,\qquad u,v\in R^{1}\nn\\
 &\le&  C_{u^p,a}(\si^p )^{2} (u-v)|\xi |,\qquad u,v\in [0,a],
\eea
for some constant $C_{u^p,a}$ depending on ${u^p}$ and $a$.

  Clearly,  
\be
  \lim_{h\to 0}\sup_{\stackrel{|u-v|\le h }{ u,v\in\De}} \frac{  (\si^p )^{2} (u-v)|\xi  |}{    (2(\si^0)^2(u-v) \log 1/ |u-v|)^{1/2} } =0,\qquad a.s., \label{5.20qq}
\ee 
where we use the fact that $(\si^{p})^{2}(x)\le (\si^{0})^{2}(x)$. Therefore, it follows from (\ref{5.37}) that,
\be
  \lim_{h\to 0}\sup_{\stackrel{|u-v|\le h }{ u,v\in\De}} \frac{ ( \ov \eta_{\bb}(u)-\ov   \eta_{\bb}(v)) a_{0} + (\ov   \eta_{p+\bb}(u)-\ov   \eta_{p+\bb}(v)) a_{p} }{    (2(\si^0)^2(u-v) \log 1/ |u-v|)^{1/2} } =1,\qquad a.s. \label{5.34}
\ee
 Using this in Theorem \ref{theo-unifmod} we get (\ref{urev.2}) in this case.

  \medskip\noindent\textbf{Case 3.}  Let  $\eta_{0}$, the Gaussian process with covariance $\frak{u}^{0}(x,y)$, and $\eta_{p}$ the Gaussian process with covariance $v^{p}(x,y)$. Then as in  Case 1. we can show that,
  \be
 \lim_{h\to 0}\sup_{\stackrel{|u-v|\le h }{ u,v\in\De}} \frac{ ( \eta_{0} (u)-\eta_{0} (v)) a_{0} + ( \eta _{ p }(u)- \eta _{p}(v)) a_{1} }{    (2(\si^0)^2(u-v) \log 1/ |u-v|)^{1/2} }=1,\quad a.s.  \label{5.27tol}
\ee
Continuing, as in    Case 2. we have,
 \be
 \lim_{h\to 0}\sup_{\stackrel{|u-v|\le h }{ u,v\in\De}} \frac{ ( \eta_{0} (u)-\eta_{0} (v)) a_{0} + (\ov \eta _{ p }(u)- \ov \eta _{p}(v)) a_{1} }{    (2(\si^0)^2(u-v) \log 1/ |u-v|)^{1/2} }=1,\quad a.s.  \label{5.27til}
\ee
Using this in Theorem \ref{theo-unifmod} we get (\ref{urev.2}) in this case.
 \qed
 \begin{example}{\rm   \textbf{Diffusions}    Let  $\mathcal{Z}$ be a  diffusion in $R^1$ that is regular and without traps  and  is symmetric with respect to a $\si$--finite measure  $m$, called the speed measure, which is absolutely continuous with respect to Lebesgue measure.  
(A diffusion is   regular and without traps when
$P^{x}\( T_{y}<\ff\)>0,   \forall x,y\in  R^1$.) We   consider diffusions $\mathcal{Z}$ with  generators of the form
\begin{equation}
\wt  L=\frac{1}{2}a^{2}(x)\frac{d^{2}}{dx^{2}}+c(x)\frac{d}{dx}\label{gendef.1},
\end{equation}
where $a^{2}(x)$ and $c(x)\in C\(R^1\)$ and $a(x)\ne 0$ for any $x\in R^1$.    For details see   \cite[Section 7.3]{RY} and \cite[Chapter 16]{Breiman}.

   Let  $s(x)$ be a scale function for $\cal Z$.  It is   strictly  increasing and  unique up to a linear transformation.  Since $a^{2}(x)$ and $c(x)\in C\(R^1\)$, $s(x)\in C^2(R^1)$.   We take $s(0)=0$  and consider the case where $\lim_{x\to \ff}s(x)=\ff$.

 For $\bb>0$, the $\bb-$potential density of $\cal Z$ with respect to  the speed measure     is,  
 \begin{equation} \label{ddiff.1}
\ov u^\bb (x,y)= \left\{
 \begin{array} {cc}
 p_\bb(x)q_\bb(y),& \quad x\leq y  
 \\
 q_\bb(x)p_\bb(y),& \quad y\leq x   
\end{array}  \right. ,
\end{equation}
where $p_\bb$ and $q_\bb$ are in $ C^{2}(R^{1})$ and  are positive and    $p_\bb$ is strictly  increasing and  
$q_\bb$ is strictly  decreasing.     See   \cite[(4.114)]{book}.  
  Let,
\begin{equation} \label{diffs.1}
  (\ov\si^\bb)^2(x,y)=\ov u^\bb (x,x)+\ov u^\bb (y,y)-2\ov u^\bb (x,v)
\end{equation}
 We show in  \cite[Lemma 4.9 and Theorem 9.1]{MRLIL} that for all $x,y\in R^1$,
\begin{equation} \label{}
  (\ov\si^\bb)^2(x,y)\sim s'(x)|x-y|,\qquad\text{as $|x-y|\to 0$}.
\end{equation}
  It follows from this that,
\begin{equation} \label{diffs.2}
  (\ov\si^\bb)^2(x,y)\sim  |s(x)-s(y)|,\qquad\text{as $|x-y|\to 0$}.
\end{equation}

  (To connect the scale function $s(x)$ to $ \ov u^\bb(x,y)$ note that for all $\bb>0$,
\begin{equation} \label{5.44}
   s'(x)=q_\bb(x)p'_\bb(x)-p_\bb(x)q'_\bb(x) .
\end{equation}
 See e.g., \cite[Theorem 9.1]{MRLIL}.)

\medskip \textbf{Case 4.} Let    $\ov Y_4=\{\ov Y_{4,t},t\ge 0 \}$ be a process  with state space $T_{4}=R^{1}$ that is obtained by killing  $\cal Z$ at the end of an independent exponential time with mean $1/\bb$. In this case    $\ov Y_4$ has $0-$potential density $\ov u^{\bb}$ and 
   $p-$potential density   $\ov u^{\bb+p}$.  
   
 It follows from (\ref{193.3a}) that
 for this process the function $w^p(x,y)$ in (\ref{int.1}) is,
\begin{equation} \label{5.50}
\ov  w_4^{p}(x,y)=\ov u^{\bb+p}(x,y)+\(\frac{1}{p}-\int_{R^1}  \ov u^{\bb+p}(x,z)\,dm(z)\)\frac{f_{4}(y)}{\|f_{4}\|_1},
\end{equation}
where \begin{equation} \label{dint.2}
  f_{4}(y)=\int_{ -\ff} ^{\ff} \ov  u^{\bb+p}(x,y)\,d\mu(x),
\end{equation}
   and the $L_1$ norm is taken with respect to $m$. When $\cal Z$ is conservative, that is, when $P^{y}\(\ze<\ff\)=0$ for all $y\in R^{1}$,  so that for all $y\in R^{1}$, we have,\begin{equation} \label{5.4mmdif}
  \int_{-\ff} ^\ff \ov   u^{p}(x,y)\,dy=\frac1{p}, \qquad \forall\, p>0.
\end{equation} 
we can simplify this as in (\ref{int.2simple}).

\medskip\textbf{Case 5.}    Let  $\ov Y_{5}=\{\ov Y_{5,t},t\ge 0 \}$, be a process with state space  $T_{5}=(0,\ff)$ that is obtained by starting  $\ov Y_4$  and then  killing it  the first time it hits 0. The $0-$potential density of $\ov Y_5$ is,
\begin{equation} \label{d5.6}
\ov  v^{\bb} :=\ov v^{\bb}(x,y)=\ov u^\bb(x,y)-\frac{\ov u^\bb(x,0 )\ov u^\bb(0, y)}{\ov u^\bb(0,0)};
\end{equation} 
see  \cite[4.165]{book}.

The  $p-$potential density of     $\ov Y_5 $   is  $\ov v^{\bb+p}$.    It follows from (\ref{193.3a}) that
 for these processes the function $w^p(x,y)$ in (\ref{int.1}) is,
\begin{equation} \label{d5.7}
\ov  w_5^{p}(x,y)=\ov v^{\bb+p}(x,y)+\(\frac{1}{p}-\int_{T}   \ov v^{\bb+p}(x,z)\,dm(z)\)\frac{f_{5}(y)}{\|f_{5}\|_1},
\end{equation}
where \begin{equation} \label{dint.3}
 f_{5}(y)=\int_{0} ^{\ff} \ov v^{\bb+p}(x,y)\,d\mu(x),
\end{equation}
 where   the $L_1$ norm is taken with respect to $m$.  

  Let,
\begin{equation} \label{diffs.1b}
  (\ov\si_{v}^\bb)^2(x,y)=\ov v^\bb (x,x)+\ov v^\bb (y,y)-2\ov v^\bb (x,v).
\end{equation}
 We show in  \cite[Lemma 4.11 and Theorem 9.1]{MRLIL} that for all $x,y\in R^1$,
\begin{equation} \label{}
 (\ov\si_{\ov v}^\bb)^2(x,y)\sim (\ov\si^\bb)^2(x,y)\sim s'(x)|x-y|,\qquad\text{as $|x-y|\to 0$}.
\end{equation}
  It follows from this that,
\begin{equation} \label{diffs.2b}
(\ov\si_{\ov v}^\bb)^2(x,y)\sim   (\ov\si^\bb)^2(x,y)\sim  |s(x)-s(y)|,\qquad\text{as $|x-y|\to 0$}.
\end{equation}

\medskip\textbf{Case 6.}  
 Let $\ov Y_{6}=\{\ov Y_{6,t},t\ge 0 \}$ 
 be a process  with state space  $T_{6}=(0,\ff)$ that is obtained by starting  $\cal Z$ in $T_{6} $ and then  killing it  the first time it hits 0.  The $0-$potential density of $\ov Y_{6}$ with respect to the speed measure  $m$
is,
  \begin{equation} \label{diff.5}
\ov{\frak{u}}^{0}(x,y)=s(x)\wedge s(y),\qquad  x,y>0.
\end{equation}
(We get this by taking the limit $a $ goes to 0 and  $b $ goes to infinity in \cite[VII. Corollary 3.8]{RY}.)

 By  \cite[Example 4.5.4]{book}, for any $p>0$ the $p-$potential density $ \ov{\frak{u}}^{p}(x,y)$ of     $\ov Y_6 $ satisfies
\begin{equation}
 \ov{\frak{u}}^{p}(x,y)= \ov v^{ p}(x,y).\label{duisv.1}
\end{equation}
 For these processes the function $w^p(x,y)$ in (\ref{int.1}) is given by  
  \be  \label{5.61}
 \ov  w_6^{p}(x,y)= \ov v^{ p}(x,y)+\(\frac{1}{p}-\int_{T}   \ov v^{ p}(x,z)\,dm(z)\)\frac{f_{6}(y)}{\|f_{6}\|_1},
\end{equation}
where \begin{equation} \label{dint.3}
 f_{6}(y)=\int_{0} ^{\ff}  \ov  v^{p}(x,y)\,d\mu(x)
\end{equation}
and the $L_1$ norm is taken with respect to $m$. 

Let  
\begin{equation} \label{diffs.1c}
  (\ov\si_{\ov{\frak{u}}^0})^2(x,y)=\ov{\frak{u}}^{0}(x,x)+\ov{\frak{u}}^{0}(y,y)-2\ov{\frak{u}}^{0}(x,y).
\end{equation}
Clearly
\begin{equation} \label{diffs.2c}
 (\ov\si_{\ov{\frak{u}}^0})^2(x,y)=  |s(x)-s(y)|,\qquad \forall\, x,y.
\end{equation}
 }\end{example}

Let $\{\wt Y_{i},i=4,5,6 \}$ denote the 
rebirthing of the process $\{  \ov  Y_{i},i=4,5,6 \}$ with probability measure $\mu$ suported on $T_{i} $. Let  
 $\wt L_{i,t}^{y}$ denote the local time of   $ \wt Y_i$
  normalized so that,
\begin{equation}
  \wt E^{ x}\(\int_{0}^{\ff}  e^{-ps}\, d_{s}\wt L^{y}_{i,s}\)=   w_i^p(x,y),\qquad x,y\in T_{i} . \label{d5.12}
\end{equation}

  \begin{theorem} \label{theo-d5.1} Consider the examples in  Cases
4., 5., and 6. Then for $i=4,5,6$, and any $d\in T_{i}$, 
   \be
   \limsup_{x\to d}\frac{|\wt L_{i,t}^{x}-\wt L_{i,t}^{d}|}{(2|s(x)-s(d)|\log\log  1/|x-d|)^{1/2}}= \(2\wt L_{i,t}^{d}\)^{1/2},  \quad a.e.\,\,\, t,  \,\,\, \wt P^{y}\,\,
a.s. \label{drev.2n}
\ee
for all $y\in T_{i}$,
and
  for any closed interval $\De$ in $T_{i}\cap [0,1]$,
   \be
\lim_{h\to 0}\sup_{\stackrel{|u-v|\le h }{ u,v\in\De}} \frac{|\wt L_{i,t}^{u}-\wt L_{i,t}^{v}|}{  (2|s(u)-s(v)| \log  1/| u-v|)^{1/2} }= \sup_{u\in \De}\(2\wt L_{i,t}^{u}\)^{1/2},  \quad a.e.\,\,\, t,  \,\,\, \wt P^{y}\,\,
a.s.\label{durev.2}
\ee
for all $y\in T_{i}$.  
\end{theorem}

\medskip \Proof  Let  $\wt G_{0}=\{\wt G_{0}(x),x\geq 0 \}$ be a mean zero Gaussian process with covariance 
 \be \label{5.68}\ov{\frak{u}}^{0}(x,y)=s(x)\wedge s(y).
 \ee
Then,   by (\ref{diffs.1c}) and (\ref{diffs.2c})
\begin{equation} \label{5.67}
E(\wt G_0(x)-\wt G_0(y))^2=|s(x)-s(y)|.
\end{equation} 
It is clear that  $\{\wt G_{0}(x),x\geq 0 \}=\{B_{s(x)},x\geq 0 \}$ where $B_t$ is Brownian motion. Therefore,   for any   $d\geq 0$,  
\be
\limsup_{x\to d}
\frac{|\wt G_{0}(x)-\wt G_{0}(d)|}{ (2(|s(x)-s(d)|)\log \log 1/|x-d|)^{1/2}}= 1\qquad \mbox {a.s.},
\label{tag3.1}
\ee
where we use,
\begin{equation} \label{}
  \lim_{x\to 0}\frac{\log 1/|x-d|}{\log 1/|s(x)-s(d)|}=1.
\end{equation}

 It follows from (\ref{diffs.2}) and (\ref{diffs.2c}) that for  $d\geq 0$  and all  $p>0$,
   \begin{equation}
 \lim_{x,y\to d } \frac{(\ov \si^p)^2(x,y)}{  \ov \si^2_{\frak {\ov u}^0}(x,y)}=1\label{1.4am}.
   \end{equation}
   
  Using (\ref{tag3.1}) and Theorem \ref{theo-localmod}, we get (\ref{drev.2n}) in Case 4.

\medskip The same argument but using  (\ref{diffs.2b}) and   (\ref{diffs.2c}) gives (\ref{drev.2n}) in Cases 5. and 6.

\medskip We now consider the uniform modulus of continuity.
 We begin with Case 6.   Let $\wt G_{0}$ be as given in (\ref{5.68}) and let  $  G_{p}=\{G_{p}(x),x\geq 0\}$  be a mean zero Gaussian process with  covariance    $\ov{\frak{u}}^{p}(x,y)$ independent of $\wt G_{0}$. By Theorem \ref{theo-unifmod}, to obtain  (\ref{durev.2}) in   Case 6 it suffices to show that   for any $a_{1}, a_{2}$ with $a^{2}_{1}+a^{2}_{2}=1$,
\be
 \lim_{h\to 0}\sup_{\stackrel{|u-v|\le h }{ u,v\in\De}} \frac{ \(\wt G_{0}(u)-\wt G_{0}(v)\) a_{1}+\(G_{p}(u)-G_{p}(v)\) a_{2} }{  (2(|s(u)-s(v)|)\log 1/|u-v|)^{1/2} } =1, \quad a.s. \label{ex.1r6x} 
\ee
  To see this note that the denominator $\vf(s(u,v))$ in Theorem \ref{theo-unifmod} is 
\bea
&&(2(|s(u)-s(v)|)\log 1/|u-v|)^{1/2}\\
&&\qquad\nn \le (2\sup_{u\in\De}s'(u)(| u-v|)\log 1/|u-v|)^{1/2}:=\wt\vf(C|u-v|),
\eea
for some constant $C$, which satisfies (\ref{3.39}).

  We now describe $G_p$. The next lemma is given in \cite[Lemma 3.3.6]{book}. We repeat it because of its significance in this proof.
\begin{lemma}
\label{lem-posdef2} Let $\{
p_t(x,y)\,;\,(t,x,y)\in
R_+\times S\times S \}$ be a family of symmetric regular transition
densities with respect to some reference measure
$m$.   (I.e.,  $p_t(x,y)$ are jointly measurable and satisfy the  Chapman--Kolomogorov identities.) Then, for each
$t> 0$,
$p_t(x,y)$ is positive definite. Hence, whenever $u^{
\al}(x,y)$ exists, for any
$\al\geq 0$, it is positive definite  and therefore is the covariance of a Gaussian process.
\end{lemma}

\Proof    Using    the Chapman--Kolmogorov equation   and
then symmetry we get
\begin{eqnarray}
\sum_{i,j=1}^n a_ia_jp_t(x_i,x_j)&=& \sum_{i,j=1}^n a_ia_j
\int p_{t/2}(x_i,z)p_{t/2}(z,x_j)\,dm(z)\nn\\ &=& \int
\big|\sum_{i=1}^n a_i p_{t/2}(x_i,z)\big|^2 \,dm(z)\geq 0.\label{mp3.6}
\end{eqnarray} This shows that $p_t(x,y)$ is positive definite. Since,
\begin{equation} \label{}
  u^{\al}(x,y)=\int_{-\ff}^{\ff}e^{-\al t}p_t(x,y)\,dt
\end{equation} 
we see that $ u^{\al}(x,y)$ is positive definite.\qed

Let $\{
 \wt p_t(x,y)\,;\,(t,x,y)\in
R_+\times T\times T \}$ be a family of symmetric regular transition
densities for $\ov Y_{6}$ with respect to speed measure. We then have that for  any $\al>0$
\begin{eqnarray}
\ov{\frak{u}}^{0}(x,y)&=&\int_{0}^{\ff} \wt p_t(x,y)\,dt
\label{tag3.6}
\\
&=&\int_{0}^{\ff}e^{-\al t} \wt p_t(x,y)\,dt+ \int_{0}^{\ff}\(1-e^{-\al t}\) \wt p_t(x,y)\,dt
\nonumber\\
&=& \ov{\frak{u}}^{\al}(x,y)+ \int_{0}^{\ff}\(1-e^{-\al t}\) \wt p_t(x,y)\,dt:=\ov{\frak{u}}^{\al}(x,y)+\ov{\frak{h}}^{\al}(x,y).
\nonumber
\end{eqnarray}
 It follows from Lemma \ref{lem-posdef2} that $\ov{\frak{h}}^{\al}(x,y)$ is also positive definite.
  Let  $H_{p}=\{H_{p}(x), x\geq 0\}$ be a mean zero Gaussian process with
covariance $\ov{\frak{h}}^{p}(x,y)$ that is independent of $G_{p}$. It follows that,

\begin{equation}
\wt G_{0}\stl G_{p}+H_{p}.\label{tag3.7}
\end{equation}
We show in  (\ref{duisv.1}) that  $\ov{\frak{u}}^{p}(x,y)= \ov v^{ p}(x,y) $.
 Using this and (\ref{5.67}) we see that,  
\be  \ov \si_{\ov{\frak{u}}}^2(x,y)=(\ov \si_{\ov v}^p ) ^{2}(x,y) +E\(H_p(x)-H_p(y)\)^2.
\ee
Consequently it follows from   (\ref{diffs.2b}) and (\ref{diffs.2c}) that,
\begin{equation} \label{5.76}
\ov \si ^2_{H_{p}}(x,y) :=E\(H_p(x)-H_p(y)\)^2=o(|x-y|),\qquad \text{as $x,y\to 0.$}
\end{equation}
Using  \cite[Corollary 7.2.3]{book} we see that for all $\ep>0$,
\be
\lim_{\de\to 0}\sup_{\stackrel{|u-v|\le\de}{u,v\in \De} }
\frac{|H_{p}(u)-H_{p}(v)|}{ (2(\ep|u-v|)\log 1/|u-v|)^{1/2}}\le 1\qquad \mbox {a.s.}
\label{5.77}
\ee
In addition since $\wt G_0$ is a time changed Brownian motion,
    for any interval $\De$ in $[0,1]$,
\be
\lim_{\de\to 0}\sup_{\stackrel{|u-v|\le\de}{u,v\in \De} }
\frac{|\wt G_{0}(u)-\wt G_{0}(v)|}{ (2(|s(u)-s(v)|)\log 1/|u-v|)^{1/2}}= 1\qquad \mbox {a.s.}
\label{tag3.1qq}
\ee
Let $\wt G'_{0}=\{\wt G'_{0}(t),t\in [0,1]\}$  
be an independent copy of $\wt G_{0}$. 
Then when $a^{2}_{1}+a^{2}_{2}=1$, it follows from (\ref{tag3.7}) that,
\bea\label{5.79}
&& \wt G_{0}(u)- \wt G_{0}(v)\nn\\
&&\qquad  \stl\(\wt G_{0}(u)-\wt G_{0}(v)\) a_{1}+ \(\wt G'_{0}(u)-\wt G'_{0}(v)\) a_{2}\\
&&\qquad \stl \(\wt G_{0}(u)-\wt G_{0,1}(v)\) a_{1}+\(G_{p}(u)-G_{p}(v)\) a_{2}+ \(H_{p}(u)-H_{p}(v)\) a_{2}.\nn 
\eea
Using (\ref{5.77}) and  (\ref{tag3.1}) we obtain (\ref{ex.1r6x}), which as we have stated,
gives (\ref{durev.2}) in Case 6. 

\medskip  The same analysis that proves Case 6. shows that when $G_{\bb }=\{G_{\bb }(t),t\in [0,1]\}$ and $G_{\bb+p }=\{G_{\bb+p }(t),t\in [0,1]\}$ are independent,  
 \be
 \lim_{h\to 0}\sup_{\stackrel{|u-v|\le h }{ u,v\in\De}} \frac{ \(  G_{\bb}(u)- G_{\bb}(v)\) a_{1}+\(G_{\bb+p}(u)-G_{\bb+p}(v)\) a_{2} }{  (2(|s(u)-s(v)|)\log 1/|u-v|)^{1/2} } =1, \quad a.s. \label{ex.1r5} 
\ee
 Using this and  we get (\ref{durev.2}) in   Case 5. 

 \medskip   Lastly we consider Case 4.  Let $\wt \eta_{\bb }=\{\wt \eta_{\bb }(t),t\in [0,1]\}$ be  a mean zero Gaussian process with covariance $\bar u^{\bb } (x,y)$,  and $\wt \eta_{\bb+p}=\{\wt \eta_{\bb+p}(t),t\in [0,1]\}$ be a mean zero Gaussian process with covariance $\bar u^{\bb+p} (x,y)$, independent of $\wt \eta_{\bb}$. To use  Theorem \ref{theo-unifmod} to get (\ref{durev.2}) in  this case, it suffices to show that, 
 \be
 \lim_{h\to 0}\sup_{\stackrel{|u-v|\le h }{ u,v\in\De}} \frac{ \(  \wt \eta_{\bb}(u)- \wt \eta_{\bb}(v)\) a_{1}+\(\wt \eta_{\bb+p}(u)-\wt \eta_{\bb+p}(v)\) a_{2} }{  (2(|s(u)-s(v)|)\log 1/|u-v|)^{1/2} } =1, \quad a.s. \label{ex.1r4} 
\ee

Let $\xi $ be an independent 
standard normal random variable.   It follows from  (\ref{d5.6}) that, 
\begin{equation}
\wt \eta_{\bb }(x)\stl G_{\bb }(x)+ c_{\bb }q _{\bb }(x) \xi ,\quad  c_{\bb}=  \({ p_\bb(0)}/{q_\bb(0)}\)^{1/2} ,\label{5try.1}
\end{equation}
as in   (\ref{5.26}) and (\ref{5.26a}). Note that since   $q_\bb$ is in $ C^{2}(R^{1})$ ,
\begin{equation} \label{}
  |q _{\bb }(u)-q _{\bb }(v)|\le C|u-v|,\qquad u,v\in [0,1],
\end{equation}
for some constant $C$. In addition, since
 $s\in C^{2}(R^1)$   and  is   strictly increasing and  ,  \begin{equation}
 |s(u)-s(v)|\ge C'    |u-v|,  \qquad u,v\in [0,1]\label{ex3.m}.
\end{equation}
Consequently, 
\be
\lim_{\de\to 0}\sup_{\stackrel{|u-v|\le\de}{u,v\in \De} }
\frac{  |q _{\bb }(u)-q _{\bb }(v)|\,|\xi|}{ (2( |s(u)-s(v)|)\log 1/|u-v|)^{1/2}}= 0,\qquad \mbox {a.s.}
\label{5.774}
\ee
 Using this and (\ref{ex.1r5}) we get (\ref{ex.1r4}).
\qed

\begin{remark}{\rm \label{rem-5.1} It follows from (\ref{5.22}) that for all $\bb\ge 0$,  $(\si^\bb)^2(x)$ is regularly varying at 0 with positive index and consequently is dominated by a regularly varying function  with positive index on $[-T,T]$ for some interval $T$. Obviously this also holds for Gaussian processes with increments variances that are asymptotic to             $(\si^\bb)^2(x)$,   such as $(\si_{v }^p)^2$ in (\ref{5.20vu}).
In \cite[Theorem 7.2.1]{book} we give upper bounds on $[-T,T]$ for the uniform modulus of continuity     of Gaussian processes with these increments variances. In particular this shows that they are continuous on $[-T,T]$. For the processes we consider we can cover   $R^1$ with these intervals and see that the Gaussian processes are continuous on  $R^1$.
}\end{remark}

 \section{Partially rebirthed Markov processes}\label{sec-partreb}

   A partially rebirthed   Markov process  is  obtained by starting a transient Markov process with state space $S$ and  when it dies returning it to $S$  with a sub-probability measure $\Xi$. (With probability    $1-|\Xi |$  it is sent to a disjoint state space $S'$, where it remains.) 
    The potential density of such a process on $S$ has the   simple form, \begin{equation}
  \wt    u^{0}(x,y)= u^{0}(x,y)+f(y), \hspace{.2 in}x,y\in S,\label{75.9jr}
\end{equation}
where $u^{0}(x,y)$ is the potential density of the original transient process and 
\begin{equation}
f(y)=\int_{S} u^{0}\(x,y\)\, d\nu(x),\label{6.2mm}
\end{equation}
 and $\nu(x)=\Xi(x)/(1-|\Xi(S)|)$.

  Let $S$ a be locally compact space with a countable base.  
Let   $\YY=(\Om,  \FF_{t},  \YY_t,\newline\th_{t},  P^x)$ be a transient symmetric Borel right process with state space $S$  and  continuous strictly positive  $p-$potential   densities  
\be
  u^{p}=\{ u^{p}(x,y), x,y\in S \},\quad p\geq 0,
\ee
 with respect to some $\si$--finite  positive measure $m$ on $S$. Let $\ze=\inf\{t\,|\,\YY_t=\De \}$,  where $\De$ is the cemetery state for   $\YY$  and assume       that $\ze<\ff$   almost surely.   

 Let $S'$ a be locally compact space with a countable base, which is disjoint from S.  
Let   $\YY'=(\Om',  \FF'_{t},  \YY'_t,\th'_{t},  P^{'x})$ be a transient symmetric Borel right process with state space $S'$,   infinite lifetime    and  continuous strictly positive  $p-$potential   densities   
\be
 ( u')^{ p}=\{   ( u')^{ p}(x,y), x,y\in S' \},\quad p\geq 0,
\ee
 with respect to some $\si$--finite  positive measure $m'$ on $S'$. (For an  example of such a process take $S'=(0,\ff)$  and $\YY' $ to be a  $\text{BES}^{3}$ process. For this process $(u')^{0}(x,y)=\min \(1/x, 1/y\) $ with respect to $m'(dy)=2y^2 dy$. 
(See \cite[p. 306]{RY}.)

\medskip  We now define a Borel right process
 $\wt X\!=\!
(\wt \Om,  \wt \FF_{t},\wt  X_t, \wt \th_{t},\wt P^x
)$ with state space $T=S\cup S' $. On $S$ the process   behaves like $\YY $, on $S'$ it   behaves like $\YY' $.  Rather than restarting  the process  at each death time,  which by assumption   only occurs when the process is in $S$,   we restart it with respect to a  probability measure $\Xi$ on $T=S\cup S' $ which we now describe.
 
 \medskip Let $\nu$  be a finite   positive measure on $S$ with mass $|\nu|=:|\nu (S)|$  and let $x'_{0}$ be a fixed point in $S'$. On $S$,     \begin{equation}
\Xi\(dw\)=\frac{\nu\(\,dw\)}{1+|\nu |}. \label{75.5}
 \end{equation}
On $S'$,  
  \begin{equation}
\Xi\( \{x'_{0}\}\)=1/(1+|\nu |), \qquad  \text{ and } \qquad   \Xi\(S'- \{x'_{0}\}\)=0. \label{75.3}
 \end{equation} 
 At each death time there is probability $1/(1+|\nu |)$ of restarting at $  x'_{0}$, in which case the process never returns to $S$.  We call the process $\wt X\!=\!
(\wt \Om,  \wt \FF_{t},\wt  X_t, \wt \th_{t},\wt P^x
)$, restricted to $S$, a partial rebirthing of $\cal Y$. 

\begin{lemma} \label{} On $S$, $\wt X\!=\!
(\wt \Om,  \wt \FF_{t},\wt  X_t, \wt \th_{t},\wt P^x
)$ has potential densities 
\begin{equation}
  \wt    u^{0}(x,y)= u^{0}(x,y)+f(y), \hspace{.2 in}x,y\in S.\label{75.9}
\end{equation}
  
\end{lemma}

\Proof Define the probability measures  $\{\nu_{i},i\ge 1\}$ on $T$ by,
\begin{equation}
\nu_{i}\(\,dw\)=\(|\nu |/(1+|\nu |)\)^{i-1}\frac{\nu\(\,dw\)}{|\nu |}, \qquad w\in S, \label{75.6a}
\end{equation}
\begin{equation}
\nu_{i}\(\{x'_{0}\}\)= 1-\(|\nu |/(1+|\nu |)\)^{i-1}\quad \text{ and } \quad \nu_{i}\(S'-\{x'_{0}\}\)=0.\label{75.6b}
\end{equation}
Then for each $i>1$ and $B\in \mathcal{M}( F(S))$
 \begin{equation}
\wt P^{y}\(B(\om_{i})\)=P^{\nu_{i}}(B)=\(|\nu |/(1+|\nu |)\)^{i-1}P^{\nu/ |\nu |}(B), \label{75.7}
 \end{equation}
where we use the notation given in the second paragraph below (\ref{2.2}). 
 Note that under $\nu_{i}$ all mass not in $S$ is placed at $x'_{0}$.

 It follows from (\ref{75.6a})--(\ref{75.7}) for any bounded measurable function $h$ on $S$,
 \begin{eqnarray}
 \lefteqn{\wt E^{x}\(\int_{0}^{\ff}h\(\wt  X_t\)\,dt \)}\\
 &&=\int_S u^{0}\(x,y\) h(y) \,dm(y)+\sum_{i=2}^{\ff}\wt E^{x}\(\int_{0}^{\ze(\om_{i}) } h(\om_{i})(t)\,dt\)
 \nonumber
 \\
 &&=\int_S u^{0}\(x,y\) h(y) \,dm(y)+\sum_{i=2}^{\ff}\(|\nu|/(1+|\nu |)\)^{i-1}E^{\nu/ |\nu |}\(\int_{0}^{\ff}h\(\YY_t\)\,dt\).\nn
\eea
This is equal to\bea   \lefteqn{\int_S u^{0}\(x,y\) h(y) \,dm(y)}\\
 &&\qquad+\sum_{i=2}^{\ff}\(|\nu|/(1+|\nu |)\)^{i-1}\frac{1}{|\nu|}\int_S       \int_{0}^{\ff} u^{0}\(x,y\) h(y) \,dm(y)\, d\nu(x)
 \nonumber \\
 &&=\int_S u^{0}\(x,y\) h(y) \,dm(y)+ \int_0^\ff     \(\int_S u^{0}\(x,y\)\, d\nu(x) \)   h(y) \,dm(y).
 \nonumber
 \end{eqnarray}
 This proves the lemma.\qed
  
  Let $L^{y}_{t}  $ be the local time for $\YY$ normalized so that,
  \begin{equation}\label{6.14}
  E^{x}\( L_{\ff}^{y}\)= u^{0}(x,y), \hspace{.2 in}x,y\in S.
  \end{equation} 
  
   \bl\label{lem-ltdecomp6}   The process $\wt X$ has a jointly continuous local time $\{\wt L_{t}^{y}, y\in S,t\in R_+ \}$ such that
for each $r\ge 1$ and $y\in S$,
  \begin{equation}
 \wt L^{y}_{t}(\wt\om)=\sum^{r-1}_{i=1}L^{y}_{\ze (\om_{i})}(\om_{i})+ L^{y}_{t-\ze_{r-1}}(\om_{r}), \quad \forall t\in(\ze_{r-1},\ze_{r}], 
 \quad a.s.\label{62.12t}
 \end{equation} 
 and for all $x,y\in S$,
 \begin{equation}
  \wt E^{ x}\(\wt L^{y}_\ff\)=   \wt u^0(x,y). \label{6az.1vw}
\end{equation}
  \el

\medskip\Proof    (\ref{62.12t}) appears to be the same statement as (\ref{2.12t}), but of course the local times are of different processes. Nevertheless, the proof of this lemma is essentially exactly the same as the proof of Lemma \ref{lem-ltdecomp}. It follows from (\ref{75.7}) that for the process restricted to $S$, $\wt P^{y}=  P_{1}^{y}\times_{n=2}^{\ff}P_n^{\nu_{i}}$. In the paragraph preceding (\ref{193.3a}) we point out that for the fully rebirthed process $\wt P^{y}=  P_{1}^{y}\times_{n=2}^{\ff}P_n^{\mu }$. However, the specific nature of this probabilty  plays no role   in the proof of Lemma \ref{lem-ltdecomp} in which we show that (\ref{2.12t}) is a jointly continuous local time. Consequently, the proof of 
Lemma \ref{lem-ltdecomp} also proves that (\ref{62.12t}) is a local time of the partially rebirthed process  and is jointly continuous.

   To see that $A_\ff^y$ has the normalization in (\ref{6az.1vw})
 note that for $x,y\in S$,
\begin{eqnarray}
  \wt E^{x}\( A_{\ff}^{y}\)&=&\sum_{i=1}^{\ff} \wt E^{x}\( L^{y}_{\ff}(\om_{i})\)
\label{caf.22}
\\
&=&E^{x}\( L^{y}_{\ff}\)+\sum_{i=2}^{\ff} \wt E^{x}\( L^{y}_{\ff}(\om_{i})\)=u^{0}(x,y)+\sum_{i=2}^{\ff} E^{\mu_{i}}\( L^{y}_{\ff}(\om_{i})\)
\nonumber\\
&=&u^{0}(x,y)+\sum_{i=2}^{\ff} \(|\nu |/(1+|\nu |)\)^{i-2}\frac{1}{1+|\nu |}\int_S  E^{z}\( L^{y}_{\ff} \)\,d\nu (z).
\nonumber\eea
This is equal to,
\be 
 u^{0}(x,y)+\sum_{i=2}^{\ff} \(|\nu |/(1+|\nu |)\)^{i-2}\frac{1}{1+|\nu |}\int_S  u^{0}(z,y)\,d\nu (z)\nn\\
=\wt  u^{0}(x,y).\nonumber
\ee
\qed

 It is easy to see that the  analogue of Lemma \ref{lem-condind} also holds for $\wt X$. 

\medskip To obtain the  analogue of Theorem \ref{theo-1.1} we have to be more careful. We want to decribe the local times $\wt L^{x}_{\ff}$ for $x\in S$. However, to employ the    isomorphism theorems we use, we must deal with probabilities. (The statement in    (\ref{3.38}) is obtained  by integrating both sides of (\ref{it1.2f}) and (\ref{it1.2})  in $y$  with respect to various probability measures.      The measures     $\nu_{i}$ in $\wt P^{y}=  P_{1}^{y}\times_{n=2}^{\ff}P_n^{\nu_{i}}$, when restricted to $S$ are   sub-probabilities.) We deal with this in the following  way:

  \medskip Let   $\YY_{+}=(\Om_{+},  \FF_{+,t},  \YY_{+,t},\th_{+,t},  P_{+}^x)$ be a transient symmetric Borel right process with state space $T=S\cup S'$, which is $\YY$ on $S$, and $\YY'$ on $S'$. $\YY_{+}$ has continuous strictly positive  $p-$potential   densities, $p\geq 0$  
\bea
  u^{p}_{+}(x,y) &=&  u^{p} (x,y), \qquad x,y\in S, \label{6.17}\\
   u^{p}_{+}(x,y') &=& 0, \hspace{.75in} x\in S, y'\in S', \nn\\
      u^{p}_{+}(x',y) &=&0, \hspace{.75in} x'\in S', y\in S , \nn\\
 u^{p}_{+}(x',y') &=&  (u')^{p} (x',y'), \hspace{.21in} x',y'\in S',  \nn
\eea
 with respect to a  positive measure $\ov m$ which is $m$ on $S$ and  $m'$ on $S'$.
 Let $L^{y}_{+,t}  $ be a local time for $\YY_{+}$ normalized so that,
  \begin{equation}
  E_{+}^{x}\( L_{+,\ff}^{y}\)= u_{+}^{0}(x,y), \quad x,y\in T.
  \end{equation}
Starting in $S$, $\YY_{+,t}$ behaves exactly like $\YY_{t}$; i.e.,  for $z\in S$, 
\be L_{+,\ff}^{z}=L_{\ff}^{z}.
\ee
 Let $\la$ denote  an independent exponential random variable with mean $1/p$
 independent of everything else. Then  
  \begin{equation}
  E_{+, \la}^{x}\( L_{+,\la}^{y}\)= u_{+}^{p}(x,y), \qquad x,y\in T.\label{6.20}
  \end{equation}
  Let 
 $ \{L^{x}_{i, t}\} _{i=1}^\ff$    be independent copies of $L^{x}_{+, t}$. Following the 
the proof of Theorem \ref{theo-1.1} exactly, we obtain:

 \bt\label{theo-6.1} Let $y\in S$. If for some $r\geq 2$ and some measurable set of functions $B\in \mathcal{M}( F(S))$ 
 \begin{equation}
\(P_{1}^{y}\times \prod^{r-1}_{i=2} P_{+, i}^{\nu_{i}}\times  P_{+, r, \la_{r} }^{\nu_{r}}\)  \( \sum^{r-1}_{i=1} L^{x}_{i, \ff}+  L^{x}_{r, \la_{r}}\in B\)=1,  \label{75.16}
 \end{equation}
 then,  
  \begin{equation}
\wt P_{\la}^{y} \(\wt L^{x}_{\la}(\wt\om)\in B\,\Big |\,    \ze_{r-1}<\la\leq \ze_{r}\)=1.\label{75.17}
 \end{equation}
 \et

\medskip  We now give an analogue of Theorem \ref{theo-ITcond} for partially rebirthed Markov processes. Let $\eta_{i, 0} (x)$, $x\in T$, i=1,2,\ldots\, be independent  Gaussian processes with covariance $ u_{+}^{0}(x,y)$ and  $\eta_{r, p} (x)$ be a Gaussian process with covariance $  u_{+}^{p}(x,y)$ independent of the $\eta_{i, 0} (x)$.   It follows from (\ref{6.17}) that  the Gaussian processes on $S$ and $S'$ are independent.  

\bt\label{theo-6.2} Set
\be 
\wt G_{r, s}\(x\)=\sum_{i=1}^{r-1} \frac{1}{2}(\eta_{i, 0 }(x)+s)^2+\frac{1}{2}(\eta_{r,p}(x)+s)^2, \label{2.21pp}
\ee
and $ P_{\wt G_{r, s}}=  \times_{i=1}^{r-1}  P_{\eta_{i, 0} } \times   P_{\eta_{r,p}} $. When
\begin{equation}
P_{\wt G_{r, s}}\(\wt G_{r, s}(\cdot)\in B\)=1,  \label{70.6na}
 \end{equation}
  for a measurable set of functions $B\in \mathcal{M}( F(S))$, then 
  \begin{equation}
\(\wt P_{\la}^{y}\times P_{\wt G_{r, s}}\) \(\wt L^{\cdot}_{\la}(\wt\om)+\wt G_{r, s}(\cdot)\in B\,\Big |\,    \ze_{r-1}<\la\leq \ze_{r}\)=1.\label{70.7n}
 \end{equation}
\et
   
  \Proof  
Set 
\be  \eta_{i, 0} (\nu_{i})=\int_T \eta_{i, 0} (x)\,d\nu_{i}(x) , \qquad  \eta_{r, p} (\nu_{r})=\int_T \eta_{r, p} (x)\,d\nu_{r}(x).
\ee  
Note that  the measures   $\nu_{i}$  are all probability measures on $T$.  This   is important and is the reason that we extended the state space $S$ to    $T$.  Recognizing this we see that        (\ref{3.38}) now takes the form:  for each  $r\geq 2$,
\bea
&&\label{4it1.3}
\Big\{ \sum^{r-1}_{i=1} L^{x}_{i, \ff}+  L^{x}_{r, \la_{r}}+\sum_{i=1}^{r-1}\frac{1}{2}(\eta_{i, 0}(x)+s)^2+\frac{1}{2}(\eta_{r, p}(x)+s)^2\,;\,\, x\in S,\,\,  \\
&&\hspace{1.5 in} \( P^{y} \times_{i=2}^{r-1}P_{+,i}^{\nu_{i}}\times  P_{+,r, \la_{r}}^{\nu_{r}} \)     \times_{i=1}^{r-1}  P_{\eta_{i, 0} } \times   P_{\eta_{r,p} } \Big\}\nn\\
 &&\stackrel{law}{=}
\Big\{ \sum_{i=1}^{r-1}\frac{1}{2}(\eta_{i, 0}(x)+s)^2+\frac{1}{2}(\eta_{r, p}(x)+s)^2\,\,\,\,\,;\,x\in S\,,\nn\\
&& \hspace{.3 in}\qquad     \,(1+\frac{\eta_{ 1}(y)}s)P_{\eta_{0,1} } \times_{i=2}^{r-1}  (1+\frac{\eta_{i, 0}(\nu_{i})}s)P_{\eta_{i, 0} }\times   ( (1+\frac{\eta_{r, p}(\nu_{r})}s)P_{\eta_{r, p} }) \Big\}.\nn
\eea
The case of $r=1$ is precisely (\ref{it1.2}). 
The rest of the proof proceeds as in the proof of Theorem \ref{theo-ITcond}.\qed

 \subsection{Moduli of continuity of the local of  times of certain partially rebirthed Markov processes}\label{sec-exampq}

 So far the material in this section gives analogues for partially rebirthed Markov processes of the material in Section \ref{sec-fullreb} for fully rebirthed Markov processes. In Section \ref{sec-app}
the results in Section \ref{sec-mod} on chi--squared processes are used in Theorems \ref{theo-localmod} and \ref{theo-unifmod}
to give local and uniform moduli of continuity for the local times of fully rebirthed Markov processes. We now note that these results hold verbatim for partially rebirthed Markov processes. To see this with respect to Theorem \ref{theo-localmod} consider (\ref{2.21}). This statement also hold with $G _{r,s}$ replaced by  $\wt G _{r,s}$ in (\ref{2.21pp}). A similar observation show that Theorem \ref{theo-unifmod} holds for partially rebirthed Markov processes. Consequently, Theorems \ref{theo-5.1} and \ref{theo-d5.1} hold as stated but with $\wt L^x_{i,t}$ the local times of the the partially rebirthed Markov processes in Cases 1.--6. There is one difference. These local times are normalized as in (\ref{6.2mm}).

To be more explicit, in Case i. the normalization is, 
\begin{equation} \label{}
  E^x(\wt L^y_{i,\ff})=\frak w_i(x,y)+\int_{T_i } \frak w_i(x,y)\,d\nu(x),
\end{equation}
where 
 \be \begin{array}{cc} 
\frak w_1(x,y)=u^\bb(x,y)\quad \text{in}\quad (\ref{1.21nn})&\quad\frak w_2(x,y)=v^\bb(x,y)  \quad \text{in}\quad(\ref{5.6})\\\\
\frak w_3(x,y)=\frak u^0(x,y) \quad \text{in}\quad (\ref{5.4})&\quad \frak w_4(x,y)=\ov u^\bb(x,y) \quad \text{in}\quad(\ref{ddiff.1})\\\\
 \frak w_5(x,y)=\ov v^\bb(x,y) \quad \text{in}\quad (\ref{d5.6})&\quad\frak w_6(x,y)= \frak{\ov u}^0(x,y)  \quad \text{in}\quad  (\ref{diff.5})
\end{array}
\ee
and 
 \be \begin{array}{ccc} 
  T_1 =R^1 &\quad T_2=R^1-\{0 \}   &\quad T_3=R^1-\{0 \} \\\\  T_4 = R^1&\quad T_5=   (0,\ff)&\quad T_6=(0,\ff).
 \end{array}
\ee

    \bibliographystyle{amsplain}

\begin{thebibliography}{99}

  \bibitem{Breiman} L.  Breiman, {\em Probability}, Classics  in Applied Mathematics, SIAM, Philadelphia,  (1992).

 
  \bibitem{Bertoin} J.  Bertoin, {\em Levy Processes}, Cambridge Tracts in Mathematics 121, Cambridge University Press, New
York,  (1996).

 \bibitem{book} M. B.  Marcus and J.~Rosen, {\em Markov Processes,
Gaussian Processes and Local Times}, Cambridge University Press, New
York,  (2006).


 \bibitem{MRECP} M. B.  Marcus and J.~Rosen, {\em  Local and uniform moduli of continuity of chi-square processes},  Electron. Commun. Probab. 27 (2022), article no. 31, 1-10. 
 
 
 
 \bibitem{MRLIL} M. B.  Marcus and J.~Rosen, {\em  Law of the iterated logarithm for $k/2$--permanental processes and the local times of    related Markov processes}, Memoirs of the AMS, to appear.
 

 


\bibitem{Meyer} P.- A. Meyer,
Renaissance, recollements, melanges, ralentissement de processes de Markov,  {\it Ann. de L'inst. Fourier},\, 25, (1975), 465--497.

 \bibitem{RY}
Daniel Revuz and Marc Yor.
\newblock {\em Continuous martingales and {B}rownian motion}, volume 293.
\newblock Springer-Verlag, Berlin, third edition, 1999.

 
 
 
\end{thebibliography}

\bigskip
\noindent
\begin{tabular}{lll} &   P.J. Fitzsimmons  \\
&Department of Mathematics  \\
& University of California, San Diego \\
& La Jolla CA, 92093 USA  \\ 
& pfitzsim@ucsd.edu\\ 
& &\\
& & \\
& Michael B. Marcus\\
&  253 West 73rd. St., Apt. 2E, \\
&  New York, NY 10023, USA \\
&mbmarcus@optonline.net\\
& &\\
& & \\
& Jay Rosen\\
& Department of Mathematics\\
&  College of Staten Island, CUNY\\
& Staten Island, NY 10314, USA \\
& jrosen30@optimum.net
\end{tabular}

\end{document}